\documentclass[a4paper, 11pt, english]{article}

\usepackage{amsmath,amsthm,amssymb,mathrsfs,mathtools,stmaryrd} 
\usepackage{graphicx}
\usepackage[left=2cm,right=2cm,top=2cm,bottom=2cm]{geometry}

\usepackage[
backend=biber,
style=alphabetic,
sorting=nyt,
backref=true,
url=false,
isbn=false,
giveninits=true,
clearlang=true]{biblatex}
\renewbibmacro{in:}{
  \ifentrytype{article}{}{\printtext{\bibstring{in}\intitlepunct}}
  }

  \DeclareFieldFormat*{title}{\mkbibemph{#1}}
  \DeclareFieldFormat*{journaltitle}{#1}
  \DeclareFieldFormat[article]{volume}{\mkbibbold{#1}}
  \DeclareFieldFormat[article]{number}{\bibstring{number}~#1}

  \renewbibmacro*{journal+issuetitle}{
    \usebibmacro{journal}
    \setunit*{\addspace}
    \iffieldundef{series}
      {}
      {\newunit
        \printfield{series}
    \setunit{\addspace}}
    \printfield{volume}
    \setunit{\addspace}
    \usebibmacro{issue+date}
    \setunit{\addcolon\space}
    \usebibmacro{issue}
    \setunit{\addcomma\space}
    \printfield{number}
    \setunit{\addcomma\space}
    \printfield{eid}
    \newunit
}

\DeclareNameAlias{author}{family-given}
\DeclareNameAlias{editor}{family-given}
\DeclareNameAlias{translator}{family-given}

\usepackage{bm}
\usepackage{enumitem}
\usepackage{tensor}
\usepackage{tikz-cd}
\usepackage{makecell}
\usepackage{hyperref}
\hypersetup{
  colorlinks=true,
  linkcolor={red!50!black},
  filecolor=magenta,
  urlcolor=cyan,
  citecolor={blue!80!black}
}
\DeclareFontFamily{U}{mathx}{\hyphenchar\font45}
\DeclareFontShape{U}{mathx}{m}{n}{
      <5> <6> <7> <8> <9> <10>
      <10.95> <12> <14.4> <17.28> <20.74> <24.88>
      mathx10
      }{}
\DeclareSymbolFont{mathx}{U}{mathx}{m}{n}
\DeclareFontSubstitution{U}{mathx}{m}{n}
\DeclareMathAccent{\wideparen}{0}{mathx}{"75}

\theoremstyle{plain}

\newtheorem{theorem}     [equation]  {Theorem}
\newtheorem{proposition} [equation]  {Proposition}
\newtheorem{lemma}       [equation]  {Lemma}
\newtheorem{corollary}   [equation]  {Corollary}

\newtheorem{innermaintheorem}        {Theorem}
\newenvironment{maintheorem}[1]
{\renewcommand\theinnermaintheorem{#1}\innermaintheorem}
{\endinnermaintheorem}

\newtheorem*{theorem*}               {Theorem}
\newtheorem*{proposition*}           {Proposition}
\newtheorem*{lemma*}                 {Lemma}
\newtheorem*{corollary*}             {Corollary}

\theoremstyle{definition}

\newtheorem{definition/} [equation]  {Definition}
\newenvironment{definition}{ \pushQED{\qed}\begin{definition/}} {\popQED\end{definition/}}

\newtheorem*{definition*/}           {Definition}
\newenvironment{definition*}{ \pushQED{\qed}\begin{definition*/}} {\popQED\end{definition*/}}

\theoremstyle{remark}

\newtheorem{remark/}     [equation]  {Remark}
\newenvironment{remark}{ \pushQED{\qed}\begin{remark/}} {\popQED\end{remark/}}
\newtheorem{example/}    [equation]  {Example}
\newenvironment{example}{ \pushQED{\qed}\begin{example/}} {\popQED\end{example/}}

\newtheorem*{remark*/}               {Remark}
\newenvironment{remark*}{ \pushQED{\qed}\begin{remark*/}} {\popQED\end{remark*/}}
\newtheorem*{example*/}              {Example}
\newenvironment{example*}{ \pushQED{\qed}\begin{example*/}} {\popQED\end{example*/}}
\newtheorem{notation/}     [equation]  {Notation}
\newenvironment{notation}{ \pushQED{\qed}\begin{notation/}} {\popQED\end{notation/}}

\newtheorem*{question*}              {Question}

\numberwithin{equation}{section}

\newenvironment{enumeratea}{
  \begin{enumerate}[label = (\alph*)]}{
  \end{enumerate}}

\newcommand{\ol}[1]{\overline{#1}}

\newcommand{\bb}[1]{\mathbb{#1}}
\newcommand{\cl}[1]{\mathcal{#1}}
\newcommand{\sr}[1]{\mathscr{#1}}
\newcommand{\fk}[1]{\mathfrak{#1}}

\newcommand{\ti}[1]{\tilde{#1}}

\newcommand{\rra}{\rightrightarrows}
\newcommand{\xar}[1]{\xrightarrow[]{#1}}
\newcommand{\fiber}[2]{\tensor[_{#1}]{\times}{_{#2}}}

\newcommand{\define}[1]{\emph{#1}}

\ExplSyntaxOn

\newcommand{\cat}[1]{
  \mathcal{\tl_range:nnn {#1} {1} {1}}~
  \text{\tl_range:nnn {#1} {2} {-1}}
  }

  \ExplSyntaxOff
  
\newcommand{\R}{\mathbb{R}}
\newcommand{\C}{\mathbb{C}}
\newcommand{\Z}{\mathbb{Z}}

\DeclareMathOperator{\colim}{colim}
\DeclareMathOperator{\dom}{dom}

\DeclareMathOperator{\ev}{ev}

\DeclareMathOperator{\GL}{GL}
\DeclareMathOperator{\id}{id}

\DeclareMathOperator{\Lie}{Lie}

\DeclareMathOperator{\pr}{pr}
\DeclareMathOperator{\per}{per}

{}

\begin{document}


\title{Lie algebras of quotient groups}
\author{David Miyamoto \footnote{Queen's University, Kingston, Canada.
    Email:
    \texttt{d.miyamoto@queensu.ca}}}


\date{\today}

\maketitle

\begin{abstract}
  We give conditions on a diffeological group $G$ and a normal subgroup $H$ under which the quotient group $G/H$ differentiates to a Lie algebra for which $\Lie(G/H) \cong \Lie(G)/\Lie(H)$. Our Lie functor is instantiated by the tangent structure on elastic diffeological spaces introduced by Blohmann. The requisite conditions on $G$ and $H$ hold, for example, when $G$ is a convenient infinite-dimensional Lie group and $H$ is countable, or when $G$ is finite-dimensional and $H$ is arbitrary. To recognize that convenient infinite-dimensional manifolds are elastic diffeological spaces, we give a characterization of convenience in terms of the diffeological tangent functor: a separated and bornological locally convex topological vector space $E$ is convenient if and only if the natural map $E \times E \to TE$ is an isomorphism of diffeological spaces. As an application, we integrate some classically non-integrable Banach-Lie algebras to diffeological groups.
\end{abstract}

\tableofcontents


\section{Introduction}
\label{sec:introduction}

One of the fundamental results in the study of Lie groups is Lie's third theorem, proved in its present form by Cartan.
\begin{theorem*}[Lie-Cartan]
  For every finite-dimensional Lie algebra $\fk{g}$, there is some Lie group whose Lie algebra is $\fk{g}$.
\end{theorem*}
Of course, not every Lie algebra is finite-dimensional, and Lie's third theorem famously fails in the infinite-dimensional case: van Est and Korthagen \cite{EstKor64} gave the first example of a Banach-Lie algebra which is not the Lie algebra of any Banach-Lie group. We call such Lie algebras $\fk{g}$ non-enlargeable, and otherwise say that $G$ integrates or enlarges $\fk{g}$. Van Est and Korthagen justify their example by identifying a topological obstruction to enlargeability within the center $\fk{z}$ of $\fk{g}$. Precisely, they showed that a Banach-Lie algebra $\fk{g}$ is enlargeable if and only if a relevant subgroup $\Pi \leq (\fk{z}, +)$ is topologically discrete. In that case, their construction of an integrating $G$ requires the Banach-Lie group structure on $\fk{z}/\Pi$, and the fact that its Lie algebra is isomorphic to $\fk{z}$.

While it is only a manifold if $\Pi$ is topologically discrete, the quotient $\fk{z}/\Pi$ is always a well-behaved diffeological group. Souriau \cite{Sour80} introduced diffeological groups, and later diffeological spaces,\footnote{K.\ T.\ Chen \cite{Chen77} gave a similar notion.} to axiomatize smoothness for non-Lie groups. The structure of a diffeological space is lightweight: a diffeology on a set $X$ is a collection of maps $U \xar{p} X$ from open subsets of Cartesian spaces, called plots, satisfying three axioms.\footnote{(i) locally constant maps are plots; (ii) if $U' \xar{F} U$ is smooth, and $U \xar{p} X$ is a plot, so is $pF$; (iii) if $p$ is locally a plot, it is a plot. One may recognize that the last two conditions describe a sheaf, and the first states that the sheaf is concrete.}  Diffeological spaces form a category $\cat{Dflg}$ which contains the category of Banach manifolds. But unlike a category of manifolds, $\cat{Dflg}$ is a quasi-topos, meaning it is complete, co-complete, and locally Cartesian closed. In particular, if $G$ is any diffeological group, and $H$ is any normal subgroup, then both $H$ and $G/H$ inherit natural diffeologies, whence $\fk{z}/\Pi$ is a diffeological group. This motivates us to integrate Banach-Lie algebras to diffeological groups. In doing so, we accomplish three objectives, which we now summarize.

\subsection{Summary of results}
\label{sec:summary-results}

First, we restrict our attention to elastic diffeological groups, introduced by Blohmann \cite{Bloh24, Bloh24b}. Whereas the existence of a Lie functor defined for all diffeological groups is unknown, we have one on the full subcategory of elastic groups:
  \begin{equation*}
    \Lie\colon \{\textrm{elastic Lie groups}\} \to \{\textrm{elastic Lie algebras}\}.
  \end{equation*}
  This Lie functor invokes a tangent functor $T\colon \cat{Dflg} \to \cat{Dflg}$, which we define as the left Kan extension of the usual tangent functor on manifolds. When restricted to the subcategory of elastic spaces, $T$ participates in a tangent structure in the sense of Rosick\'{y} \cite{Ros84}, and more precisely a Cartesian tangent structure with scalar $\R$-multiplication in the sense of \cite{AinBloh25}. It is the tangent structure which imbues $\Lie(G) \coloneqq T_eG$ with its Lie bracket. We show that, under some conditions, the Lie functor is compatible with quotients.
  \begin{maintheorem}{I}
    Let $G$ be an elastic group, and let $H$ be a normal subgroup of $G$. If both $H$ and $G/H$ are elastic, and if iterated tangents $T^k\iota\colon T^kH \to T^kG$ of the inclusion $\iota\colon H \to G$ are inductions for all $k \geq 0$, then the natural map $\Lie(G)/\Lie(H) \to \Lie(G/H)$ is an isomorphism of Lie algebras.
  \end{maintheorem}
The theorem applies, for instance, when $H$ is diffeologically discrete (e.g.\ countable), or if $G$ is any finite-dimensional Lie group. To verify that $G/H$ is an elastic group in the former case, we use the auxiliary result that to be elastic is preserved under the presence of a ``Q-chart,'' whose name we inherited from Barre \cite{Bar73}. More precisely, we establish in Theorem \ref{thm:groupoids-2} that if $\varphi \colon X \to Y$ is a Q-chart, and $X$ is elastic, then so is $Y$. In lieu of defining Q-charts here, we give the motivating example: if $H$ is a countable group acting freely and smoothly on a diffeological space, then the quotient map $X \to X/H$ is a Q-chart. For a concrete example, $\R \to \R/(\Z + \alpha\Z)$ is a Q-chart when $\alpha$ is irrational.

Second, we embed convenient Lie groups, in the sense of Kriegl and Michor \cite{KriegMic97}, into the tangent category of elastic diffeological groups. These infinite-dimensional Lie groups are modelled by convenient locally convex topological vector spaces, and the category of convenient manifolds includes the categories of Banach and Fr\'{e}chet manifolds. We have found that convenience is intimately related to properties of the tangent bundle. Precisely, we call a diffeological vector space $V$ \define{tangent-stable} if the map
\begin{equation*}
  V \times V \to TV, \quad (u,v) \mapsto \frac{d}{dt}\Big|_{t=0}u+tv
\end{equation*}
is an isomorphism. We show:
\begin{maintheorem}{II}
  A separated bornological locally convex topological vector space is convenient if and only if, when viewed as a diffeological space, it is tangent-stable.
\end{maintheorem}

Finally, we treat enlargeability of convenient Lie algebras. We follow the strategy used by Cartan, and extended by Neeb \cite{Neeb02}, of integrating central extensions of Lie algebras. Given a topologically split central extension $\fk{a} \hookrightarrow \hat{\fk{g}} \twoheadrightarrow \fk{g}$ of convenient Lie algebras, we may find a corresponding smooth cocycle $\omega \in Z_s^2(\fk{g}, \fk{a})$. If $G$ is a connected and simply-connected convenient Lie group with $\Lie(G) = \fk{g}$, we get a period homomorphism $\per_\omega \colon H_2(G) \to \fk{a}$. We prove:

\begin{maintheorem}{III}
  Suppose that $\fk{a} \hookrightarrow \hat{\fk{g}} \twoheadrightarrow \fk{g}$ is a topologically split central extension of convenient Lie algebras. Fix $\omega \in Z_s^2(\fk{g}, \fk{a})$ corresponding to this extension. Assume that there is some connected and simply-connected convenient Lie group $G$ with $\Lie(G) = \fk{g}$. If the image $\Pi_\omega$ of the associated period homomorphism $\per_\omega\colon H_2(G) \to \fk{a}$ is diffeologically discrete, then there is some elastic diffeological group $\hat{G}$ centrally extending $G$ by $\fk{a}/\Pi_\omega$, for which
  \begin{equation*}
    \begin{tikzcd}
      \fk{a}/\Pi_\omega \ar[d, dashed, "\Lie"] \ar[r, hook] & \hat{G} \ar[r, two heads] \ar[d, dashed, "\Lie"] & G \ar[d, dashed, "\Lie"]  \\
      \fk{a} \ar[r, hook] & \hat{\fk{g}} \ar[r, two heads] & \fk{g}.
    \end{tikzcd}
  \end{equation*}
\end{maintheorem}
As an application, we integrate some non-enlargeable Lie algebras, such as the classical examples given by van Est and Korthagen \cite{EstKor64} and Douady and Lazard \cite{DouadLaz66}.

\subsection{Remarks on the literature}
\label{sec:remarks-literature}

There is a litany of diffeological tangent functors. They broadly fall into three approaches: the internal, representing tangent vectors by curves with equality probed by representations in higher dimensional plots; the external, identifying tangent vectors with derivations of rings of functions; and the mixed, representing tangent vectors by curves, but testing equality using smooth functions. We have traced some treatments of the various functors below.

\begin{table}[h]
  \centering
  \begin{tabular}[h]{c | c | c }
    Internal & External & Mixed \\ \hline 
    Losik \cite{Los92} & Iglesias-Zemmour \cite{Igl13} & Souriau \cite{Sour80} \\
    Hector \cite{Hec95} & Christensen, Wu \cite{ChrWu16} & Torre \cite{Tor98}  \\
    Christensen, Wu \cite{ChrWu16}  & Taho \cite{Tah24} &  Leslie \cite{Les03} \\
            Blohmann \cite{Bloh24} & & Dugmore, Ntumba \cite{DugNtum07} \\
             & & Vincent \cite{Vin08} \\
             & & Stacey \cite{Stac13} \\
    & & Taho \cite{Tah25}
  \end{tabular}
  \caption{Some approaches to tangent functors}
  \label{tab:1}
\end{table}

We also mention Magnot \cite{Mag18}, who follows Dugmore and Ntumba, and Goldammer and Welker \cite{GolWel21}, who follow Vincent. We use Blohmann's tangent functor, which coincides with Losik's. Its closest, yet distinct, relative is Christensen and Wu's. Like us, Leslie and Magnot assign Lie algebras to some diffeological groups. However, it is unclear whether our desired examples satisfy their hypotheses.

Tangent structures were introduced by Rosick\'{y} \cite{Ros84}, and Cockett and Cruttwell (e.g.\ \cite{CocCrut14,CocCrut15}) have substantially expanded the theory. Recently, Aintablian and Blohmann \cite{AinBloh25} established a Lie functor that operates on groupoid objects in tangent categories, and that is the Lie functor we use on the category of elastic diffeological groups.

The standard reference for convenient analysis is Kriegl and Michor's book \cite{KriegMic97}, though the categorical treatment of Fr\"{o}licher and Kriegl \cite{FroelKrieg88} more directly reveals the connection to diffeology, via Fr\"{o}licher spaces. Hain \cite{Hain79} showed that the category of Banach manifolds embeds into $\cat{Dflg}$. Losik \cite{Los92} showed the same for Fr\'{e}chet manifolds, and also showed compatibility of the tangent functors. The embedding of convenient manifolds into $\cat{Dflg}$ is known to experts, and treated explicitly in e.g.\ \cite[Chapter 2]{Kih23}. However, the fact that the tangent functor for convenient manifolds coincides with the diffeological tangent functor is new, as is the characterization of convenience in terms of tangent-stability.

Methods for integrating central extensions of Lie algebras originate with Cartan \cite{Car52}, with advances by van Est and Korthagen \cite{EstKor64} and Neeb \cite{Neeb02}. Our integration into diffeological spaces is new, but the principal of using more flexible categories to avoid obstructions is not. Wockel, individually \cite{Woc11} and later with Zhu \cite{WocZhu16}, integrated central extensions of locally exponential locally convex Lie algebras into central extensions of certain Lie 2-groups. This works without caveats, at the expense of higher categorical structures. Belti\c{t}\u{a} and Pelletier \cite{BelPel25}, following an approach based on integrating paths instantiated by {\'S}wierczkowski \cite{Swier71}, integrate Banach-Lie algebras into what they call H-manifolds, a generalization of Banach Q-manifolds (which themselves are spaces locally arising as the codomain of a Q-chart from a Banach space). These are all sets with structure, in contrast with the approach via higher structures. In earlier work, Plaisant \cite{Plais80} used the same approach to integrate separable Banach Lie algebras into Banach Q-manifolds, although without the framework of a Lie functor. 
\subsection{Structure of the paper}
\label{sec:structure-paper}

In Section \ref{sec:tangent-structure}, we introduce the tangent category of elastic spaces, with a particular focus on elastic groups and vector spaces. In Section \ref{sec:diff-group}, we address the elasticity of orbit spaces of diffeological groupoids, and in particular establish that the codomains of Q-charts from elastic spaces are also elastic. We then apply our general result to subgroups of elastic groups, and achieve Theorem \ref{thm:groupoids-1}. In Section \ref{sec:conv-sett-diff}, we introduce the convenient setting of analysis from a diffeological perspective, and prove Theorem \ref{thm:convenient-1}. Finally, in Section \ref{sec:integr-centr-extens}, we integrate some central extensions of convenient Lie algebras, and in particular integrate specific Banach-Lie algebras to diffeological groups via Theorem \ref{thm:central-extensions-1}. In Appendix \ref{sec:tangent-categories} we review tangent categories.

\subsection*{Acknowledgements}

I am grateful to Christian Blohmann for introducing me to tangent structures and elastic spaces. I would also like to thank Lory Aintablian and David Aretz for helpful discussions on tangent categories and convenient analysis, respectively. Finally, I would like to acknowledge the Max Planck Institute for Mathematics in Bonn, which provided a friendly and supportive environment while I pursued this work.

\section{A tangent structure}
\label{sec:tangent-structure}

This section has two parts. In Subsection \ref{sec:elast-diff-spac}, we introduce our tangent functor $T\colon \cat{Dflg} \to \cat{Dflg}$, and the tangent category of elastic spaces $\cat{Elst}$. In Subsection \ref{sec:elast-diff-groups-1}, we restrict our focus to diffeological vector spaces. We emphasize that it is the entire tangent structure on $\cat{Elst}$ that lets us define the Lie functor on elastic groups.

We assume basic familiarity with diffeological spaces: the definition, distinguished maps such as inductions and subductions (in categorical language, the strong monomorphisms and strong epimorphisms), and the D-topology. For the necessary background, the reader may refer to the first two chapters of Iglesias-Zemmour's book \cite{Igl13}. For a sheaf-based introduction to diffeology, we recommend Baez and Hoffnung's article \cite{BaezHof11}.

\subsection{Elastic diffeological spaces}
\label{sec:elast-diff-spac}

Let $\cat{Cart}$ denote the category of open subsets of Cartesian (Euclidean) spaces, with smooth maps between them. We have the full and faithful functor $y\colon \cat{Cart} \to \cat{Dflg}$, which assigns to an open subset of Cartesian space its canonical diffeology.

\begin{definition}
  Given a functor $\hat{F}\colon \cat{Cart} \to \cat{Cart}$, we denote the left Kan extension of $y\hat{F}$ along $y$ by $\bb{L}\hat{F} \coloneqq \operatorname{Lan}_y y\hat{F}$.
\end{definition}

We can also take the left Kan extension of a natural transformation. Given a natural transformation $\hat{\alpha}\colon \hat{F} \to \hat{G}$, we denote the left Kan extension of $y\hat{\alpha}\colon y\hat{F} \to y\hat{G}$ along $y$ by $\bb{L}\hat{\alpha} \coloneqq \operatorname{Lan}_y y\hat{\alpha}$. This is a natural transformation from $\bb{L}\hat{F}$ to $\bb{L}\hat{G}$. We will now describe how to compute $(\bb{L}\hat{F})_X$ and $(\bb{L}\hat{\alpha})_X$.

The category $\cat{Cart}$ is small and $\cat{Dflg}$ is co-complete. Therefore, for any functor $\hat{F}$, the left Kan extension $\bb{L}\hat{F}$ exists, and moreover it is given by the colimit
  \begin{equation*}
    (\bb{L}\hat{F})X = \colim_{y \downarrow X} y\hat{F} \Pi^X.
  \end{equation*}  
  Here we are taking the colimit over the category of plots of $X$, which is the comma category $y \downarrow X$. Explicitly, $y \downarrow X$ has
  \begin{itemize}
  \item objects: pairs $(U,p)$, where $p\colon yU \to X$ is a plot;
  \item morphisms: smooth maps $f \colon U \to U'$ that fit into the diagram
    \begin{equation*}
       \begin{tikzcd}
        yU & & yU' \\
        & X. &
        \ar["yf", from=1-1, to=1-3]
        \ar["p"', from=1-1, to=2-2]
        \ar["p'", from=1-3, to=2-2]
      \end{tikzcd}
    \end{equation*}
  \end{itemize}
The functor $\Pi^X\colon y \downarrow X \to \cat{Cart}$ is the projection sending $(U,p)$ to $U$.

In order to compute $(\bb{L}\hat{F})X$, we review how to compute arbitrary small colimits in $\cat{Dflg}$. Given any functor $D\colon \cat{J} \to \cat{Dflg}$ from a small category $\cat{J}$, its colimit is induced by the co-equalizer
  \begin{equation*}
    \begin{tikzcd}
      \coprod_{f \in \operatorname{Mor}(\cl{J})} D(\dom f) & \coprod_{U \in \operatorname{Ob}(\cl{J})} DU  & \colim_{\cl{J}} D.
      \ar["\id", shift left, from=1-1, to=1-2]
      \ar["\coprod Df"', shift right, from=1-1, to=1-2]
      \ar["\pr", two heads, from=1-2, to=1-3]
    \end{tikzcd}
  \end{equation*}
  In detail, the two coproducts have their standard diffeology, and $\colim_{\cl{J}} D$ is the diffeological quotient of $\coprod DU$ by the relation generated by $x \sim Df(x)$. Since $Df$ is not necessarily invertible, $x \sim y$ if and only if there is a \define{zig-zag}, also called a \define{chain}, between them: a collection of elements $(x_i \in DU_i)_{i=0}^{k+1}$ with $x_0 = x$ and $x_{k+1} = y$, and morphisms $(f_i)_{i=0}^k$, where for each $0\leq i\leq k$, either $Df_i(x_i) = x_{i+1}$ or $Df_i(x_{i+1}) = x_i$.

  \begin{notation}
    \label{not:tangent-structure-1}
    We observe that $(\bb{L}\hat{F})X$ is a quotient of $\coprod_{(U,p)} \hat{F}U$. Given $u \in \hat{F}U$ indexed by $p$, we write $p_*u$ for the corresponding element in $(\bb{L}\hat{F})X$. For a natural transformation $\hat{\alpha} \colon \hat{F} \to \hat{G}$, we have $(\bb{L}\hat{\alpha})_X(p_*u) = p_* (\hat{\alpha}_U(u))$
  \end{notation}
  


We now turn to the tangent functor on $\cat{Dflg}$, and its role in a tangent category. We summarize the theory of tangent categories in Appendix \ref{sec:tangent-categories}. The Cartesian tangent structure on $\cat{Cart}$ with scalar $\R$ multiplication (see Definition \ref{def:appendix-1}) is given by the following data: the tangent functor $\hat{T}\colon \cat{Cart} \to \cat{Cart}$, defined by
\begin{equation*}
  \hat{T}(U \subseteq \R^n) \coloneqq U \times \R^n, \quad \hat{T}f(u,v) \coloneqq (f(u), D_uf(v)),
\end{equation*}
where $D_uf$ is the Jacobian of $f$ at $u$, the usual ring structure on $\R$, and the natural transformations
\makeatletter
\@fleqntrue
\setlength\@mathmargin{\parindent}
 \begin{alignat*}{3}
  &\textrm{the footprint} \quad & \hat{\pi}&\colon \hat{T} \to 1 \quad & (u,v) &\mapsto u \\
  &\textrm{the zero section}\quad & \hat{0}&\colon 1 \to \hat{T} \quad & u &\mapsto (u,0) \\
  &\textrm{the scalar multiplication}\quad & \hat{\cdot} &\colon \R \times \hat{T} \to \hat{T}\quad & (t,u,v) &\mapsto (u,tv)\\
  &\textrm{the addition}\quad & \hat{+} &\colon \hat{T} \fiber{}{\hat{\pi}}\hat{T} \to \hat{T} \quad & (u,v_0,u,v_1) &\mapsto (u, v_0 + v_1) \\
  &\textrm{the symmetric structure}\quad & \hat{\tau} &\colon \hat{T}^2 \to \hat{T}^2 \quad & (u,v_0,v_1,v_{01}) &\mapsto (u,v_1,v_0,v_{01}) \\
  &\textrm{the vertical lift}\quad & \hat{\lambda}&\colon \hat{T} \to \hat{T}^2 \quad & (u,v) &\mapsto (u,0,0,v).
\end{alignat*}
\@fleqnfalse
\makeatother

Formally, we say the tuple $(\hat{T}, \hat{\pi}, \hat{0}, \hat{\cdot}, \hat{+}, \hat{\tau}, \hat{\lambda})$ is a Cartesian tangent structure on $\cat{Cart}$ with scalar $\R$ multiplication. To get such a tangent structure on $\cat{Dflg}$, we would like to apply $\bb{L}$ to this tuple.

\begin{definition}
  The \define{tangent functor} $T\colon \cat{Dflg} \to \cat{Dflg}$ is defined by $T \coloneqq \bb{L}\hat{T} = \operatorname{Lan}_yy\hat{T}$.
\end{definition}
We use this functor extensively.
\begin{notation}
  \label{not:tangent-structure-2}
We unpack Notation \ref{not:tangent-structure-1} for the tangent functor. The tangent bundle $TX$ is a quotient of
  \begin{equation*}
    \coprod_{(U,p)} y\hat{T}U = \coprod_{(U,p)} yU \times y\R^n.
  \end{equation*}
  For a generic element $(u,v)$ in this coproduct, indexed by $(U,p)$, we write $p_*(u,v)$ for the corresponding element in $TX$. We may also use the shorthand $p_*v$, and note that $p_*v = Tp(u,v)$. We may say ``$p_*v$ is the tangent vector represented by $v$ in the plot $p$.''

  To illustrate a chain witnessing an equality $(p_b)_*v_b = (p_e)_*v_e$ (using $b$ for ``beginning'' and $e$ for ``end''), we will write $(yU,v) \xar{f} (yU',v')$ to denote that $\hat{T}f(u,v) = (u',v')$. Then such a chain consists of the commutative diagram:
  \begin{equation*}
    \begin{tikzcd}
      (yU_b,v_b) & (yU_1,v_1)  & (yU_2, v_2) & \cdots & (yU_e, v_e)  \\
      & & X. & &
      \ar["f_b", from=1-1, to=1-2]
      \ar["f_1"', from=1-3, to=1-2]
      \ar["f_2", from=1-3, to=1-4]
      \ar["f_e"', from=1-5, to=1-4]
      \ar["p_b"', from=1-1, to=2-3]
      \ar["p_1", from=1-2, to=2-3]
      \ar["p_2", from=1-3, to=2-3]
      \ar["p_e", from=1-5, to=2-3]
    \end{tikzcd}
  \end{equation*}
 Any arrow may be reversed. For a smooth map $f\colon X \to Y$, its tangent $Tf\colon TX \to TY$ is given by $Tf(p_*v) = (f \circ p)_*(v)$.
\end{notation}

The tangent functor shares the following properties with the tangent functor on manifolds.
\begin{definition}
  We use $\partial$ to denote the map
  \begin{equation*}
    \partial\colon C^\infty(y\R, X) \to TX, \quad \gamma \mapsto \gamma_*(0,1).
  \end{equation*}
  We may write $\delta = \partial_X$ to indicate the space $X$, or $\delta(\gamma) = \partial_t\gamma(t)$ to include a curve's parameter.
\end{definition}

\begin{proposition}[{}]
  \label{prop:tangent-structure-1}
  The tangent functor has the following properties:
  \begin{itemize}
\item \cite[Corollary 2.2.19]{Bloh24b} It preserves products, meaning that $T(X \times Y)$ is naturally isomorphic to $TX \times TY$.
\item \cite[Proposition 2.2.16]{Bloh24b} It preserves subductions, meaning that if $\varphi\colon X \to Y$ is a subduction, so is $T\varphi\colon TX \to TY$.
\item \cite[Proposition 3.35]{Bloh24} It is local, meaning that if $U \subseteq X$ is D-open, with inclusion $i\colon U \hookrightarrow X$, then $Ti\colon TU \to TX$ is an induction with D-open image. In other words, we identify $TU$ with an open subset of $TX$.
  \item \cite[Proposition 2.2.20]{Bloh24b} The map $\partial$ is a subduction onto $TX$.
\end{itemize}
\end{proposition}

We move on to the natural transformations in the tangent structure.

\begin{definition}
  In $\cat{Dflg}$, we define the footprint projection $\pi$, the zero section $0$, and the multiplication $\cdot$, respectively as the left Kan extensions of $\hat{\pi}$, $\hat{0}$, and $\hat{\cdot}$.
\end{definition}

The footprint of $p_*(u,u_0) \in TX$ is $p(u)$. The zero section satisfies $0(x) = p_*(u,0)$, for any plot $yU \xar{p} X$ with $p(u) = x$. Finally, the scalar multiplication is $t \cdot p_*(u,u_0) = p_*(u,tu_0)$.

We have the $k$-fold fiber products $T_k \coloneqq T \times_\pi \cdots \times_\pi T$. For a tangent structure, $T_2$ must be the source of addition, and $T^2$ must be the source and target of the symmetric structure. This is where we encounter problems, because in general $\bb{L}\hat{T}_2$ is not isomorphic to $T_2$. We will see a concrete example of the first difficulty in Example \ref{ex:groupoids-1}. Furthermore, $\bb{L}(\hat{T}^2)$ is not generally isomorphic to $T^2$.

Nevertheless, we always have the natural transformations (using Notation \ref{not:tangent-structure-1})
\begin{alignat*}{2}
  \theta_k&\colon \bb{L}\hat{T}_k \to T_k \quad & &\theta_{k,X}(p_*(u,v_0,\ldots, v_{k-1})) = (p_*v_0,\ldots, p_*v_{k-1}) \\
  \theta^2&\colon \bb{L}\hat{T}^2 \to T^2 \quad & &\theta^2_X(p_*(u,v_0,v_1,v_{01})) = T^2p(u,v_0,v_1,v_{01}) \\
 \nu_k&\colon TT_k \to T^2 \times_T \cdots \times_T T^2 \quad & &\nu_{k,X}((p_1,\ldots, p_k)_*v) = ((p_1)_*v, \ldots, (p_k)_*v).
\end{alignat*}
To define a tangent structure, Blohmann \cite{Bloh24} has isolated a subcategory of $\cat{Dflg}$ on which these natural transformations behave properly.

\begin{definition}[{\cite[Definition 4.1]{Bloh24}}]
\label{def:tangent-structure-1}
  A diffeological space $X$ is \define{elastic} if it satisfies:

\begin{itemize}
\item[(E1)] The maps $\theta_{k,X}$ are isomorphisms for all $k \geq 1$.
\item[(E2)] There is map $\tau_X \colon T^2X \to T^2X$ (the \define{canonical flip}) such that
  \begin{equation*}
    \begin{tikzcd}
      (\bb{L}\hat{T}^2)X \ar[r, "(\bb{L}\hat{\tau})_X"] \ar[d, "\theta_X^2"] & (\bb{L}\hat{T}^2)X \ar[d, "\theta_X^2"] \\
      T^2X \ar[r, "\tau"] & T^2X
    \end{tikzcd}
  \end{equation*}
  commutes.
\item[(E3)] The \define{vertical lift}
  \begin{equation*}
    \begin{tikzcd}
    \lambda_X\colon TX \ar[r, "(\bb{L}\hat{\lambda})_X"] & (\bb{L}\hat{T}^2)X \ar[r, "\theta_X^2"] & T^2X
  \end{tikzcd}
\end{equation*}
is an induction.
\item[(E4)] The maps $\nu_{k,X}$ are injective for all $k \geq 1$.
  \item[(E5)] For natural numbers $k_1,\ldots, k_n$, the spaces $T_{k_1}\cdots T_{k_n}X$ satisfy (E1) through (E4).
\end{itemize}
\end{definition}

A diffeological space $X$ is either elastic or not. The collection of elastic diffeological spaces, with smooth maps between them, is a full sub-category of $\cat{Dflg}$, which we denote by $\cat{Elst}$. The \emph{raison d'\^{e}tre} for these axioms is the following theorem: 
\begin{theorem}[{\cite[Theorem 4.2]{Bloh24}}]
  \label{thm:tangent-structure-2}
  On $\cat{Elst}$, the tangent functor $T$, the ring $\R$, and the natural transformations
  \begin{equation*}
    \pi, \quad 0, \quad \cdot,\quad + \coloneqq \bb{L}\hat{+} \circ \theta_2^{-1}, \quad \tau \text{ from (E2)}, \quad \lambda,
  \end{equation*}
form a Cartesian tangent category with scalar $\R$ multiplication.
\end{theorem}

\begin{remark}
  In the next section, Example \ref{ex:groupoids-1} will show that (E1) may fail. Similarly, (E4) can fail. We also suspect that (E3) can fail, with a potential example coming from a non-convenient locally convex topological vector space (cf.\ Section \ref{sec:conv-sett-diff}). However, it may be that (E2) holds for any diffeological space. 
\end{remark}

When $G$ is an elastic group, meaning a diffeological group that is a member of $\cat{Elst}$, the tangent structure on $\cat{Elst}$ identifies $\fk{g} \coloneqq T_eG$ with $\fk{X}(G)^L$, the space of left-invariant vector fields on $G$. The tangent structure also gives $\fk{X}(G)^L$ the Lie bracket inherited from $\fk{X}(G)$, and thus we get a Lie algebra structure on $\fk{g}$. Moreover, this assignment is functorial; given a smooth homomorphism $f\colon G \to H$, its tangent map $Tf\colon \fk{g} \to \fk{h}$ is a homomorphism of Lie algebras. All of these statements depend only on the fact that $\cat{Elst}$ is a tangent category, and are consequences of the work by Aintablian and Blohmann \cite{AinBloh25} which we review in Appendix \ref{sec:tangent-categories}. The next proposition summarizes the situation. 
\begin{proposition}
  \label{prop:tangent-structure-4}
We have a Lie functor
  \begin{equation*}
    \operatorname{Lie}\colon \{\text{elastic groups}\} \to \{\text{elastic Lie algebras}\},
  \end{equation*}
  given by $\Lie(G) \coloneqq \fk{g} \cong \fk{X}(G)^L$.
\end{proposition}

The Lie bracket on $\fk{X}(G)^L$ derives from the bracket on $\fk{X}(G)$ furnished by the tangent structure on $\cat{Elst}$, as described in Appendix \ref{sec:conv-sett-diff}. However, there is a more familiar formulation of the bracket. Identifying $T(X \times Y) \cong TX \times TY$, for a function $f\colon X \times Y \to Z$, we have the partial tangent morphisms
\begin{align*}
  &
    \begin{tikzcd}[ampersand replacement = \&]
      T_{(1)}f\colon TX \times Y \ar[r, "1 \times 0"] \& TX \times TY \ar[r, "Tf"] \& TZ 
    \end{tikzcd} \\
  &    \begin{tikzcd}[ampersand replacement = \&]
    T_{(2)}f\colon X \times TY \ar[r, "0 \times 1"] \& TX \times TY \ar[r, "Tf"] \& TZ.
  \end{tikzcd} \\  
\end{align*}
Denote the conjugation on $G$ by
\begin{equation*}
  c\colon G \times G \to G, \quad (g,h) \mapsto ghg^{-1}.
\end{equation*}
\begin{lemma}
  \label{lem:tangent-structure-1}
  For an elastic Lie group $G$, and two elements $v,w\in T_eG$, we have
  \begin{equation*}
    \lambda_2(w,[v,w]) = T_{(1)}T_{(2)}c(v,w).
  \end{equation*}
  where $\lambda_2(a,b) = \tau(T0(a)+\lambda(b))$. In particular, $T_{(1)}T_{(2)}c(v,w)$ completely determines $[v,w]$.
\end{lemma}
\begin{proof}
  This holds in a general tangent category, so we prove it in Lemma \ref{lem:appendix-2}.
\end{proof}

\subsection{Elastic diffeological vector spaces}
\label{sec:elast-diff-groups-1}
In this subsection, we apply the previous theory to diffeological vector spaces. This machinery will help in Section \ref{sec:conv-sett-diff}.

\begin{definition}
  A diffeological vector space over $\R$ is a vector space object over $\R$ in $\cat{Dflg}$. In other words, a diffeological vector space is a diffeological space $V$ equipped with a vector space structure for which addition $V \times V \to V$ and scalar multiplication $\R \times V \to V$ are smooth. 
\end{definition}

The key observation is that, given any diffeological vector space $V$, the map $\partial\colon C^\infty(\R,V) \to TV$ pulls back to a map $\ell = \ell_V$ from $V \times V$, using lines as distinguished curves:
\begin{equation*}
\ell = \ell_V \colon V \times V \to TV, \quad (u,v) \mapsto \partial[u+tv] = T[u+tv](0,1).
\end{equation*}
Here we use square brackets to denote the curve $t \mapsto u+tv$, and may write $\partial = \partial_t$ to emphasize that $t$ is the parameter for the curve. We have the expected identities for the projection, zero section, and scalar multiplication:
\begin{align*}
  \pi_V(\ell(u,v)) = u, \quad 0_V(u) = \ell(u,0), \quad t\cdot \ell(u,v) = \ell(u,tv).
\end{align*}

Under the identification of $TV \times TV$ with $T(V \times V)$, we have the map $\ell^2$ defined as the composition
  \begin{equation*}
    \begin{tikzcd}[sep=small]
      \ell^2\colon (V \times V) \times (V \times V) & T(V \times V) & T^2V\\
      ((u,v_0),(v_1,v_{01})) & \partial[(u,v_0) + t(v_1,v_{01})] & \partial_t\partial_s[u+tv_1 + sv_0 + tsv_{01}]
      \ar["\ell_{V\times V}", from=1-1, to=1-2]
      \ar["T\ell_V", from=1-2, to=1-3]
      \ar[mapsto, from=2-1, to=2-2]
      \ar[mapsto, from=2-2, to=2-3]
    \end{tikzcd}
  \end{equation*}
For $\vec{v} = (u,v_0,v_1,v_{01})$, we define the parametrized surface $S^2(\vec{v})(s,t) \coloneqq u+sv_0+tv_1+stv_{01}$. This surface satisfies
      \begin{equation*}
        \ell^2(\vec{v}) = T^2S^2(\vec{v})(0,e_1,e_2,0) = \theta^2(S^2(\vec{v})_*(0,e_1,e_2,0)),
      \end{equation*}
      where $e_i$ are the standard basis vectors of $\R^2$. The map
      \begin{equation*}
        N^2\colon V^4 \to (\bb{L}\hat{T}^2)V, \quad \vec{v} \mapsto S(\vec{v})_*(0,e_1,e_2,0)
      \end{equation*}
      is smooth. Therefore $\ell^2$ lifts along $\theta^2$ via $N^2$.     
  
      We also have the map $\ell_k$ defined by
      \begin{equation*}
    \ell_k\colon V^{k+1} \to T_kV, \quad (v,v_0,\ldots, v_{k-1}) \mapsto (\ell(v,v_0),\ldots, \ell(v,v_{k-1})).
  \end{equation*}
Setting $\vec{v} = (u,v_0,\ldots, v_{k-1})$, and defining $S_2(\vec{v})(t_0,\ldots, t_{k-1})\coloneqq  u + \sum t_i v_i$, we have
      \begin{equation*}
        \ell_k(\vec{v}) = T_kS_2(\vec{v})(0,e_1,\ldots, e_k) = \theta_k(S_2(\vec{v})_*(0,e_1,\ldots, e_k)),
      \end{equation*}
      and the map
      \begin{equation*}
        N_k\colon V^{k+1} \to (\bb{L}\hat{T}_k)V, \quad \vec{v} \mapsto S_2(\vec{v})_*(0,e_1,\ldots, e_k)
      \end{equation*}
      is smooth. Therefore $\ell_k$ lifts along $\theta_k$ via $N_k$.
  \begin{lemma}
    \label{lem:tangent-structure-6}
    Let $V$ be a diffeological vector space. The following four diagrams commute.
    \begin{alignat*}{2}
      & \begin{tikzcd}[ampersand replacement =\&]
        {(\bb{L}\hat{T}_2)V} \& T_2V \\
        TV \& \\
        V \times V \& V \times V \times V
        \ar["\theta_2", from=1-1, to=1-2]
        \ar["{\bb{L}\hat{+}}"', from=1-1, to=2-1]
        \ar["\ell", from=3-1, to=2-1]
        \ar["\ell_2", from=3-2, to=1-2]
        \ar["\mathrm{sum}", from=3-2, to=3-1]
        \ar["N_2", from=3-2, to=1-1]
      \end{tikzcd} \quad
     &  &\begin{tikzcd}[ampersand replacement =\&]
          {\bb{L}(\hat{T}^2)V} \& {\bb{L}(\hat{T}^2)V} \\
          V^4 \& V^4
          \ar["{\bb{L}\hat{\tau}}", from=1-1, to=1-2]
          \ar["N^2", from=2-1, to=1-1]
          \ar["N^2", from=2-2, to=1-2]
          \ar["\mathrm{flip}", from=2-1, to=2-2]
        \end{tikzcd}  \\
        & \begin{tikzcd}[ampersand replacement =\&]
          TV \& {(\bb{L}\hat{T}^2)V} \& T^2V \\
          V \times V \& \& V^4
          \ar["\bb{L}\hat{\lambda}", from=1-1, to=1-2]
          \ar["\theta^2", from=1-2, to=1-3]
          \ar["\ell", from=2-1, to=1-1]
          \ar["\ell^2", from=2-3, to=1-3]
          \ar["\mathrm{lift}", from=2-1, to=2-3]
          \ar["N^2", from=2-3, to=1-2]
        \end{tikzcd} \quad
      & & \begin{tikzcd}[ampersand replacement =\&]
        TT_2V \& T^2V \fiber{}{TV} T^2V \\
        (V^3)^2 \& V^4 \fiber{}{\pr_{12}} V^4 
        \ar["\nu_2", from=1-1, to=1-2]
        \ar["T\ell_2 \circ \ell_{V^3}"', from=2-1, to=1-1]
        \ar["\ell^2 \times \ell^2", from=2-2, to=1-2]
        \ar["\mathrm{inj}", from=2-1, to=2-2]
      \end{tikzcd}
    \end{alignat*}
    Here
\begin{alignat*}{2}
  &\mathrm{sum}(u,v_0,v_1) = (u,v_0+v_1) \quad & &\mathrm{flip}(u,v_0,v_1,v_{01}) = (u,v_1,v_0,v_{01}) \\
  &\mathrm{lift}(u,v_0) = (u,0,0,v_0) \quad & &\mathrm{inj}(u,v_0,v_1,u',v_0',v_1') = ((u,u',v_0,v_0'), (u,u',v_1,v_1'))
\end{alignat*}

\end{lemma}
\begin{proof}
  We treat each diagram separately, but advise the reader that the arguments are all similar. For the first diagram,
    \begin{align*}
      \bb{L}\hat{+} \circ N_2(u,v_0,v_1) &= \bb{L}\hat{+} \circ S_2(\vec{v})_*(0,e_1,e_2) &&\textrm{by definition of }N_2\\
                                         &= TS_2(\vec{v})(0,e_1+e_2) &&\textrm{by definition of }\bb{L}\hat{+}\\
                                         &= \partial_t[S_2(\vec{v})(t,t)] &&\textrm{representing }e_1+e_2 \textrm{ by a curve}\\
                                         &= \partial_t[u+tv_0+tv_1] \\
                                         &= \ell(u,v_0+v_1) &&\textrm{by definition of }\nu.
    \end{align*}

    For the second diagram,
    \begin{align*}
      \bb{L}\hat{\tau} \circ N^2(u,v_0,v_1,v_{01}) &= \bb{L}\hat{\tau} \circ S^2(\vec{v})_*(0,e_1,e_2,0) &&\textrm{by definition of }N^2 \\
                                                   &= S^2(\vec{v})_*(0,e_2,e_1,0) &&\textrm{by definition of }\bb{L}\hat{\tau}\\
                                                   &= N^2(u,v_1,v_0,v_{01}) &&\textrm{by definition of }N^2
    \end{align*}

    For the third diagram,
    \begin{align*}
      \bb{L}\hat{\lambda} \circ \ell(u,v) &= [u+tv]_*(0,0,0,1) &&\textrm{by definition of }\bb{L}\hat{\lambda} \\
                                           &= S^2(u,0,0,v)_*(0,e_1,e_2,0)
    \end{align*}
    The last equality follows from the fact that $S^2(u,0,0,v)$ is the composition of the curve $[u+tv]$ with the map $(s,t) \mapsto st$, and the second tangent of the map $(s,t) \mapsto st$ sends $(0,e_1,e_2,0)$ to $(0,0,0,1)$.

    For the fourth diagram,
    \begin{align*}
      T\ell_2 \circ \ell_{V^3}(\vec{v}, \vec{v'}) &= \partial_t \ell_2(\vec{v} + t\vec{v'}) \\
                                                  &= \partial_t \partial_s [(u+tu'+s(v_0+tv_0'), u+tu'+s(v_1+tv_1'))] \\
                                                  &= (\nu^2(u,u',v_0,v_0'), \nu^2(u,u',v_1,v_1').
    \end{align*}
\end{proof}

Using the above lemma, we can say more in the case that $\ell$ is an isomorphism.
\begin{definition}
  A diffeological vector space $V$ is \define{tangent-stable} if $\ell \colon V \times V \to TV$ is an isomorphism of diffeological spaces.
\end{definition}

\begin{lemma}
  \label{lem:tangent-structure-7}
A tangent-stable vector space $V$ satisfies (E3) and (E4). If $N_k$ is surjective for all $k$, then $V$ satisfies (E1). If $N^2$ is surjective, then $V$ further satisfies (E2) with $\tau \coloneqq \ell \circ \mathrm{flip} \circ \ell^{-1}$. In the case that $V$ satisfies (E1) through (E4), it is elastic, and $\ell$ is a map of bundles of vector spaces over $V$.
\end{lemma}
\begin{proof}
  Assume that $\ell$ is an isomorphism. Then so are $\ell^2$ and $\ell_k$. Therefore the diagram for $\bb{L}\hat{\lambda}$ in Lemma \ref{lem:tangent-structure-6}, and the fact that $\mathrm{lift}$ is an induction, implies that $\lambda = \bb{L} \hat{\lambda} \circ \theta^2$ is an induction, i.e., (E3) holds. Similarly, the diagram for $\nu_2$ implies that $\nu_2$ is an injection, i.e., (E4) holds.

For the statements about (E1) and (E2), observe that $N_k \circ (\ell_k)^{-1}$ is a section of $\theta_k$, and $N^2 \circ (\ell^2)^{-1}$ is a section of $\theta^2$. Therefore if $N_k$ is onto, then $\theta_k$ has a smooth inverse, satisfying (E1), and if $N^2$ is onto, then $\theta^2$ has a smooth inverse and the proposed $\tau$ satisfies (E2).
  
  Assuming (E1) through (E4) hold, the diagram for $\bb{L}\hat{+}$, and the definition of $+\colon T_2V \to TV$, together imply that $\ell$ is an isomorphism of bundles of vector spaces over $V$. Axiom (E5) holds as a consequence of the fact that if two diffeological spaces $X$ and $Y$ satisfy (E1) through (E4), then so does their product.
\end{proof}

\section{Diffeological groupoids}
\label{sec:diff-group}

In this section we justify one of our key examples of elastic groups: certain quotients of elastic groups. This is a three-step process. In Subsection \ref{sec:tang-bundl-orbit}, we compute the tangent bundle of the orbit space of some diffeological groupoids. In Subsection \ref{sec:elast-cert-quot}, we show that if $\varphi \colon X \to Y$ is a Q-chart, and $X$ is elastic, then so is $Y$. Finally, in Subsection \ref{sec:elast-quot-groups}, we illustrate when certain quotients of elastic groups remain elastic.

\subsection{Tangent bundles of orbit spaces}
\label{sec:tang-bundl-orbit}

\begin{definition}
  \label{def:3}
  A \define{diffeological groupoid} is a groupoid object in $\cat{Dflg}$. In other words, a diffeological groupoid consists of a groupoid $\cl{G} \rra X$, and diffeological structures on $\cl{G}$ and $X$ for which all structure maps are smooth.
\end{definition}

We will denote the source map by $\alpha$, the target map by $\beta$, the multiplication by $m$, the unit by $u$, and the inversion by $i$. As a consequence of the algebraic structure of $\cl{G} \rra X$, the source and target are subductions, the inversion is a diffeomorphism, and the unit is an induction. We denote the orbit map by
\begin{equation*}
  \pr = \pr^0\colon X \to X/\cl{G}.
\end{equation*}
This is a subduction, and specifically co-equalizes the pair $(\alpha,\beta)$.
\begin{remark}
  Lie groupoids, by contrast, are more than simply groupoid objects in the category of manifolds. For a Lie groupoid, we additionally require that: the base space is Hausdorff and second-countable; the source, hence target, maps are submersions; and the source, hence target, fibers are Hausdorff. Submersions between manifolds are subductions, but the converse does not hold.
\end{remark}
Let $\cl{G}_{\bullet}$ denote the nerve of a groupoid $\cl{G}$. Thus, at level $k$, we have
\begin{equation*}
  \cl{G}_k \coloneqq \cl{G} \fiber{\alpha}{\beta} \cdots \fiber{\alpha}{\beta} \cl{G}.
\end{equation*}
We also denote $\cl{G}_1 \coloneqq \cl{G}$ and $\cl{G}_0 \coloneqq X$. We adopt the next definition from \cite[Definition 3.5]{AinBloh25}. 
\begin{definition}
  A diffeological groupoid $\cl{G}$ is \define{differentiable} if the natural morphism
  \begin{equation*}
    T^n(\cl{G}_k) \to T^n\cl{G} \fiber{T^n\alpha}{T^n\beta} \cdots \fiber{T^n\alpha}{T^n\beta} T^n\cl{G}
  \end{equation*}
  is an isomorphism for each $k \geq 0$ and $n \geq 0$.
\end{definition}
When $\cl{G}$ is differentiable, $T\cl{G}$ can be viewed as a diffeological groupoid $T\cl{G} \rra TX$, with source $T\alpha$, target $T\beta$, and multiplication $Tm$.

\begin{remark}
  Strictly speaking, Aintablian and Blohmann define differentiable groupoids in a general tangent category, in which case some more properties are required. In that context, the differentiable groupoids are precisely those which differentiate to, in their sense, Lie algebroids. 
\end{remark}

\begin{example}
  Aintablian and Blohmann showed \cite[Proposition 8.1]{AinBloh25} that a groupoid object in the category of manifolds is differentiable if and only if it is a Lie groupoid. This may be interpreted as a justification for the definition of a Lie groupoid being stronger than merely a groupoid object in the category of manifolds.
\end{example}

Action groupoids are one example of differentiable diffeological groupoids.

\begin{lemma}
  \label{lem:7}
  Let $X \rtimes H \rra X$ be an action groupoid for a diffeological group $H$ acting on a diffeological space $X$. Then $X \rtimes H$ is differentiable.
\end{lemma}
\begin{proof}
  Recall that $X \rtimes H$ has for its source and target
  \begin{align*}
    \textrm{the action map } \mu\colon X \times H \to X,& \quad (x,h) \mapsto x \cdot h \\
    \textrm{the projection map} \pr_1\colon X \times H \to X,& \quad (x,h) \mapsto x.
  \end{align*}
  Thus, for instance,
  \begin{equation*}
    (X \rtimes H)_2 \cong X \times H \times H, \quad (x,h,x' = x \cdot h,h') \mapsto (x,h,h').
  \end{equation*}
  Therefore differentiability of $X \rtimes H$ follows from the fact that $T$ preserves finite products.
\end{proof}

\'{E}tale diffeological groupoids are another example of differentiable groupoids.
\begin{definition}
  A diffeological groupoid $\cl{G} \rra X$ is \define{\'{e}tale} if its source (hence target) is a local diffeomorphism. In other words if, about every $g \in \cl{G}$, there is an open neighbourhood $A \subseteq \cl{G}$ and $B \subseteq X$ such that $s\colon A \to B$ is a well-defined diffeomorphism.
\end{definition}

\begin{lemma}
  \label{lem:groupoids-1}
\'{E}tale groupoids are differentiable.
\end{lemma}

\begin{proof}
  We treat $T\cl{G}_2$, and the general case of $T^n\cl{G}_k$ follows by induction. First, since $\alpha$ and $\beta$ in the pullback diagram
  \begin{equation*}
    \begin{tikzcd}
      \cl{G}_2 & \cl{G}  \\
      \cl{G} & X
      \ar["\pr_2", from=1-1, to=1-2]
      \ar["\pr_1"', from=1-1, to=2-1]
      \ar["\alpha", from=2-1, to=2-2]
      \ar["\beta", from=1-2, to=2-2]
    \end{tikzcd}
  \end{equation*}
  are local diffeomorphisms, so too are $\pr_1$ and $\pr_2$. Indeed, pullbacks preserve local diffeomorphisms in $\cat{Dflg}$. For the same reason, and also because $T$ preserves local diffeomorphisms, it follows that every map in the diagram below, except \emph{a priori} $\nu$, is a local diffeomorphism.
  \begin{equation*}
    \begin{tikzcd}
      T \cl{G}_2  & & \\
      & T\cl{G} \fiber{T\alpha}{T\beta} T\cl{G} & T\cl{G}  \\
      & T\cl{G} & TX
      \ar["T\pr_2", bend left, from=1-1, to=2-3]
      \ar["T\pr_1"', bend right, from=1-1, to=3-2]
      \ar["\nu", from=1-1, to=2-2]
      \ar[from=2-2, to=2-3]
      \ar[from=2-2, to=3-2]
      \ar["T\alpha", from=3-2, to=3-3]
      \ar["T\beta", from=2-3, to=3-3]
    \end{tikzcd}
  \end{equation*}
  But then \emph{a fortiori} $\nu$ is also a local diffeomorphism. Furthermore, every map in the diagram above is a map of bundles that is a bijection on the fibers, except \emph{a priori} $\nu$. But then \emph{a fortiori} $\nu$ is also a bijection on the fibers. Since $\nu$ is a bundle map over the identity that is both a local diffeomorphism and a bijection on fibers, it is a diffeomorphism. 
\end{proof}

When $\cl{G}$ is a differentiable groupoid, the orbit space $TX/T\cl{G}$ of $T\cl{G}$ is a co-equalizer of the pair $(T\alpha,T\beta)$. On the other hand, the map $T\pr\colon TX \to T(X/\cl{G})$ also satisfies $T\pr \circ T\alpha = T\pr \circ T\beta$. Thus the universal property of co-equalizers gives us a unique smooth map $\eta_{\cl{G}}$ fitting in
\begin{equation*}
  \begin{tikzcd}
    & TX  & \\
    TX/T\cl{G}  & & T(X/\cl{G}).
    \ar["\pr^1"', from=1-2, to=2-1]
    \ar["T\pr", from=1-2, to=2-3]
    \ar["\eta_{\cl{G}}", from=2-1, to=2-3]
  \end{tikzcd}
\end{equation*}

The map $\eta_{\cl{G}}$ descends to the identity on $X/\cl{G}$. We do not anticipate that this map is always an isomorphism. However, it will be an isomorphism in the cases relevant to us.

\begin{proposition}
  \label{prop:groupoids-5}
  Suppose that $\cl{G} \rra X$ is a differentiable groupoid satisfying the following condition:
  \begin{equation}
    \label{eq:6}
    \textrm{If }(q,q')\colon U \to X \fiber{}{\pr} X \textrm{ is a plot, then } (Tq,Tq')\colon TU \to TX \fiber{}{\pr^1} TX \textrm{ is a plot}.
  \end{equation}
  Then the natural map $\eta_{\cl{G}}\colon T(X/\cl{G}) \to TX/T\cl{G}$ is an isomorphism.
\end{proposition}
\begin{proof}
Since the quotient $\pr$ is a subduction, so is $T\pr$. The quotient map $\pr^1$ is also a subduction. Therefore $\eta_{\cl{G}}$ is a subduction, and it is an isomorphism if it is injective. To show injectivity, we show that if $T\pr(w_b) = T\pr(w_e)$ ($b$ for ``beginning'' and $e$ for ``end''), then $\pr^1(w_b) = \pr^1(w_e)$.

Write $w_b = (q_b)_*v_b$ and $w_e = (q_e)_*v_e$. An arbitrary link in the chain witnessing the equality $T\pr(w_b) = T\pr(w_e)$ has the form (using Notation \ref{not:tangent-structure-2})
  \begin{equation*}
        \begin{tikzcd}
      (U,u,v) & & (U',u,v') \\
      & X/\cl{G}, &
      \ar["f", from=1-1, to=1-3]
      \ar["p"', from=1-1, to=2-2]
      \ar["p'", from=1-3, to=2-2]
    \end{tikzcd}
  \end{equation*}
  Since $\pr\colon X \to X/\cl{G}$ is a subduction, we may lift $p$ and $p'$ along $\pr$ at $u$ and $u'$ to maps $q$ and $q'$ into $X$. This may involve shrinking $U$ or $U'$ about $u$ or $u'$. Then we have a plot
  \begin{equation*}
    (q, q'f)\colon U \to X \fiber{}{\pr} X 
  \end{equation*}
  which, by our assumption, implies that we have a plot
  \begin{equation*}
    (Tq, T(q'f))\colon TU \to TX \fiber{}{\pr^1} TX.
  \end{equation*}
  In particular, $\pr^1(q_*v) = \pr^1(q'_*v')$.

  At the leftmost arm of a chain witnessing $T\pr(w_b) = T\pr(w_e)$, choose $q_b$ as the lift of $\pr\circ  q_b$. At the rightmost arm, choose $q_e$ as the lift of $\pr \circ q_e$. By choosing lifts in each intermediate link of the chain as we did above, we conclude that the quantity $\pr^1((q_b)_*v_b)$ is preserved across the chain, and must therefore coincide with $\pr^1((q_b)_*v_b)$. In other words, $\pr^1(w_b) = \pr^1(w_e)$, and so $\eta_{\cl{G}}$ is injective.
\end{proof}

Orbifold groupoids satisfy the preconditions of Proposition \ref{prop:groupoids-5}. Because orbifolds are not our primary focus, we discuss them in an example.

\begin{example}
\label{ex:groupoids-1}
Let $\Z_2$ act on $\R$ by the involution $x \mapsto -x$. Suppose that $q,q'\colon U \to \R$ are two plots, such that for each $u \in U$, we have $q(u) = \pm q'(u)$. Set
\begin{equation*}
  \Delta_+ \coloneqq \{u \in U \mid q(u) = q'(u)\} \textrm{ and } \Delta_- \coloneqq \{u \in U \mid q(u) = -q'(u)\}.
\end{equation*}
Then $U = \Delta_+ \cup \Delta_-$ is a finite union of closed sets, and as a consequence of the Baire category theorem, the union $\Delta_+^\circ \cup \Delta_-^\circ$ is dense in $U$. On each $\Delta_\pm^\circ$, we have $Tq = \pm Tp$. Thus, by continuity, for each $(u,v) \in TU|_{\ol{\Delta_\pm^\circ}}$ we have $Tq(u,v) = \pm Tp(u,v)$. Since $U = \ol{\Delta_+^\circ} \cup \ol{\Delta_-^\circ}$, we have showed that condition \eqref{eq:6} holds. Therefore by Proposition \ref{prop:groupoids-5}, $T(\R/\Z_2) \cong (T\R)/\Z_2$, where $\Z_2$ acts on $T\R$ by the involution $(x,v) \mapsto (-x,-v)$.

Using a variant of Proposition \ref{prop:groupoids-5} for $\bb{L}\hat{T}_2$, we find that $\bb{L}\hat{T}_2(\R/\Z_2) \cong (T_2\R)/\Z_2$. The spaces $(T\R)/\Z_2$ and $(T_2\R)/\Z_2$ are not isomorphic, so we have a concrete example where (E1) fails. In particular, $\R/\Z_2$ is not elastic.

More generally, one can show that if $\cl{G} \rra M$ is an orbifold groupoid, meaning that it is locally isomorphic to action groupoids of the form $\Gamma \ltimes \R^n \rra \R^n$ for finite groups $\Gamma \leq \GL(n)$,\footnote{The reader may be more comfortable defining an orbifold groupoid as an \'{e}tale and proper Lie groupoid. Our definition is strictly more general, and allows for non-Hausdorff orbifolds.} then $T(M/\cl{G}) \cong TM/T\cl{G}$. Upon embedding the category of orbifolds into $\cat{Dflg}$ (as done in \cite{IglKarZad10}), we observe that our tangent functor coincides with the usual orbifold tangent functor. If any of the isotropy groups $\cl{G}_x$ are non-trivial, then $M/\cl{G}$ is not elastic.
\end{example}

We turn to the special case that $H$ is a subgroup of a diffeological group $G$, equipped with its subspace diffeology. Let $\iota\colon H \hookrightarrow G$ denote the inclusion. Let $\mu\colon G \times G \to G$ denote the multiplication in $G$. Right multiplication by $H$ gives the action groupoid $G \rtimes H \rra G$. This is differentiable, so we have its tangent groupoid $T(G \rtimes H) \rra G$. The source map is $T(\mu \circ (1 \times \iota)) = T\mu \circ (1 \times T\iota)$, and the target map is $T\pr_1$. 

The tangent bundle $TG$ of $G$ is also a diffeological group, with multiplication given by $T\mu$. Therefore, as groupoids,
\begin{equation*}
  T(G \rtimes H) \cong TG \rtimes T\iota(TH),
\end{equation*}
where $T\iota(TH)$ acts on $TG$ from the right.

\begin{corollary}
  \label{cor:1}
  Let $\iota\colon H \hookrightarrow G$ be a subgroup of a diffeological group. For the associated action groupoid $G \rtimes H \rra G$, the natural map
  \begin{equation*}
    \eta_{G \rtimes H} \colon TG/T\iota(TH) =  TG/T(G \rtimes H) \to  T(G/H) 
  \end{equation*}
  is an isomorphism.
\end{corollary}
\begin{proof}
  Write the action of $H$ on $G$ by $gh \coloneqq \mu(g,h)$. Being an action groupoid, $G \rtimes H$ is differentiable. We verify it satisfies condition \ref{eq:6}, whereby the corollary holds by Proposition \ref{prop:groupoids-5}. Let
  \begin{equation*}
    (1, d)\colon G \fiber{}{\pr} G \to G \times H \textrm{ denote the inverse of } (1, \mu)\colon G \times H \to G \fiber{}{\pr} G.
  \end{equation*}
  Here $d$ stands for ``division.'' Suppose that $(q,q')\colon U \to G \fiber{}{\pr} G$ is a plot. Recall that $\pr^1$ is the orbit map for the groupoid $T(G \rtimes H)$, which has source $T\mu$ and target $T \pr_1$. Observe that
  \begin{align*}
    &T\mu \circ (Tq, Td \circ T(q,q')) = T(\mu \circ (q, d \circ (q,q'))) = Tq' \\
    &T\pr_1 \circ (Tq, Td \circ T(q,q')) = T(\pr^1 \circ (q,d \circ (q,q'))) = Tq.
  \end{align*}
  Thus $(Tq,Tq')\colon TU \to TG \fiber{}{\pr^1} TG$ is a plot, as required.
\end{proof}
\begin{remark}
  The conclusion of Corollary \ref{cor:1} generalizes to the case that $\cl{G} = X \rtimes H \rra X$, where we assume that the anchor map $(1, \mu)$ is an isomorphism. In this case, the orbit map $\pr\colon X \to X/H$ is a \define{principal fibration with structure group }$H$ in the sense of \cite{Igl13}, and the conclusion writes that $TX/TH \cong T(X/H)$. The proof is identical.
\end{remark}

\subsection{Elasticity of certain quotients}
\label{sec:elast-cert-quot}

We now restrict our attention to a special class of quotient relations on a space $X$, for which elasticity of $X$ implies elasticity of $X/{\sim}$.

\begin{definition}
  A map $\varphi \colon X \to Y$ is a \define{Q-chart} if
  \begin{itemize}
  \item[(Q1)] it is a subduction;
  \item[(Q2)] the diagonal $X \fiber{}{1} X$ is D-open inside the fiber product $X \fiber{}{\varphi} X$;
  \item[(Q3)] the second (equivalently first) projection $X \fiber{}{\varphi} X \to X$ is a local diffeomorphism.
  \end{itemize}
\end{definition}

\begin{remark}
  The notion of a Q-chart first appeared in the work of Barre \cite{Bar73}, who required their domains to be finite-dimensional manifolds. Such Q-charts also appeared in previous work of the author \cite{Miy24b}. Plaisant \cite{Plais80} generalized Barre's Q-manifolds to allow their domains to be Banach manifolds. To handle more general quotients of Banach manifolds, Belti\c{t}\u{a} and Pelletier \cite{BelPel25} introduced the related notion of an H-chart.
\end{remark}

For any map $\varphi\colon X \to Y$, we will denote $\cl{R}_\varphi \coloneqq X \fiber{}{\varphi} X$. This has the structure of a groupoid over $X$, where the source $\alpha$ and target $\beta$ are the second and first projections of $\cl{R}_\varphi$ to $X$, and the multiplication is $(z,y)\cdot (y,x) \coloneqq (z,x)$. We call $\cl{R}_\varphi \rra X$ the \define{relation groupoid} associated to $\varphi$. If $\varphi$ is a Q-chart, then $\cl{R}_\varphi$ is an \'{e}tale groupoid, hence is differentiable. We have two main examples of Q-charts.
\begin{proposition}
  \label{prop:groupoids-2}
  If $X$ is a diffeological space, and $G$ is a countable group acting smoothly and freely on $X$, then the quotient map $\varphi\colon X \to X/G$ is a Q-chart.
\end{proposition}
\begin{proof}
  Quotient maps are subductions by definition. To see that the diagonal is D-open, suppose that $(p',p) \colon U \to X \fiber{}{\varphi} X$ is a plot, and that it sends $r_0$ into the diagonal. Take a relatively compact and connected neighbourhood $V$ of $r$. The sets
  \begin{equation*}
    \Delta_g \coloneqq \{r \in \ol{V} \mid p(r) = g \cdot p'(r)\}
  \end{equation*}
  partition $\ol{V}$ into countably many disjoint closed sets. By Sierpinski's theorem \cite{Sier18}, such a partition has at most one component. This component must be $\Delta_e$, since $r_0 \in \Delta_e$. Therefore $(p',p)$ sends $V$ into $X \fiber{}{1} X$, and so the diagonal is indeed D-open in $X \fiber{}{\varphi} X$.

  To see that the second projection $X \fiber{}{\varphi} X \to X$ is a local diffeomorphism, observe that for each $g$, its restriction to $\{(g\cdot x', x') \mid x' \in X\}$ is a diffeomorphism. This set is open because it is the image of the diagonal under the diffeomorphism of $\cl{R}_\varphi$ given by $(x,x') \mapsto (g\cdot x, x')$.
\end{proof}

\begin{proposition}
  \label{prop:groupoids-3}
  If $G$ is a diffeological group, and $H$ is a subgroup whose subdiffeology is discrete, then the quotient map $\varphi\colon G \to G/H$ is a Q-chart. 
\end{proposition}
\begin{proof}
  We check that the diagonal is D-open, so suppose that $(p,p') \colon U \to X \fiber{}{\varphi} X$ is a plot that sends $r_0$ to the diagonal. We can use the smooth division map $(1,d)$ to get a plot of $H$:
  \begin{equation*}
    \begin{tikzcd}
      U & G \fiber{}{\varphi} G & X \times H & H.
      \ar["{(p,p')}", from=1-1, to=1-2]
      \ar["{(1,d)}", from=1-2, to=1-3]
      \ar["\pr_2", from=1-3, to=1-4]
    \end{tikzcd}
  \end{equation*}
  Since $H$ has discrete diffeology, this plot must be locally constant. Its value at $r_0$ is $e \in H$, so $(p,p')$ must take values in the diagonal on a neighbourhood of $r_0$. In other words, the diagonal is D-open.

  The map $\varphi$ is a subduction, and the second projection $G \fiber{}{\varphi} G \to G$ is a local diffeomorphism, for the same reasons as in the previous Proposition \ref{prop:groupoids-2}.
\end{proof}

A \define{transition} on a diffeological space $X$, denoted $h\colon X \dashrightarrow X$, is given by D-open subsets $A$ and $B$ of $X$, and a map $h\colon A \to B$, such that $h$ is a diffeomorphism. Equivalently, a transition is a locally-defined embedding from $X$ to itself.

\begin{definition}
  Given any map $\varphi\colon X \to Y$, we set
  \begin{equation*}
    \Psi(\varphi) \coloneqq \{h\colon X \dashrightarrow X \mid \varphi \circ h = \varphi \textrm{ on } \dom h\}.
  \end{equation*}
  This is the set of transitions that preserve the fibers of $\varphi$.
\end{definition}
\begin{lemma}
  \label{lem:groupoids-3}
  Let $\varphi \colon X \to Y$ be a Q-chart. For each $(x,x') \in \cl{R}_\varphi$, there is a transition $h \in \Psi(\varphi)$ such that $h(x') = x$. Moreover, the germ of such a transition at $x'$ is unique.
\end{lemma}

\begin{proof}
  Take a Q-chart $\varphi$ and suppose that $(x,x') \in \cl{R}_\varphi$. Since the second projection $\cl{R}_\varphi \to X$ is a local diffeomorphism, take a local inverse $\sigma$ mapping $x'$ to $(x,x')$. Then the desired transition is $h \coloneqq \beta \circ \sigma$.

  Suppose that $h$ and $\ti{h}$ are two transitions in $\Psi(\varphi)$ sending $x'$ to $x$. Then the juxtaposition $(\ti{h},h)$ takes values in $\cl{R}_\varphi$, and sends $x'$ to the element $(x,x)$ of the diagonal. Since the diagonal is D-open in $\cl{R}_\varphi$, we conclude that $(\ti{h},h)$ sends a neighbourhood of $x'$ into the diagonal, or in other words that $\ti{h} = h$ in a neighbourhood of $x'$. 
\end{proof}

In fact, $\Psi(\varphi)$ is a pseudogroup of transitions of $X$, and the lemma above states that the orbits of $\Psi(\varphi)$ are precisely the fibers of $\varphi$. 

\begin{lemma}
  \label{lem:groupoids-2}
  Let $\varphi\colon X \to Y$ be a Q-chart. Then $T\varphi$ is fiberwise a bijection. In other words, for all $x \in X$, the restriction $T_x\varphi \colon T_xX \to T_{\varphi(x)}Y$ is a bijection.
\end{lemma}

\begin{proof}
  Say that $\varphi(x) = y$. For surjectivity of $T_x\varphi$, we begin with the fact that $T\varphi\colon TX \to TY$ is onto, because $\varphi$ is a subduction. Thus given $w \in T_yY$, we may find $v' \in T_{x'}X$ with $T\varphi(v') = w$. Since $\varphi(x) = \varphi(x')$, by Lemma \ref{lem:groupoids-3} we may choose a locally-defined embedding $h\colon X \dashrightarrow X$ that preserves the $\varphi$-fibers and sends $x$ to $x'$. Then $v \coloneqq Th^{-1}(v')$ is the required element of $T_xX$ with $T\varphi(v) = w$.

  For injectivity of $T_x\varphi$, suppose that $v_b,v_e \in T_xX$ and $T\varphi(v_b) = T\varphi(v_e)$. An arbitrary link in the chain witnessing this equality has the following form, with the dashed arrows to be justified below. We write the link in terms of germs.
    \begin{equation*}
    \begin{tikzcd}
      (U,r) \ar[rr, "f"] \ar[dr, dashed, "\tilde{q}"] \ar[ddr, "q"'] & & (U',r') \ar[dl, dashed, "\tilde{q'}"'] \ar[ddl, "q'"] \\
      & (X,x) \ar[d, "\varphi"] & \\
      & (Y,y). & 
    \end{tikzcd}
  \end{equation*}
  The maps $\ti{q}$ is a lift of $q$ along $\varphi$ sending $r$ to $x$. Some lift of $q$ along $\varphi$ defined near $r$ exists because $\varphi$ is a subduction, and we may choose a lift mapping $r$ to $x$ by Lemma \ref{lem:groupoids-3}. We define $\ti{q'}$ similarly. Then, the juxtaposition $(\ti{q},\ti{q'}\circ f)$ takes values in $\cl{R}_\varphi$, and maps $r$ to the element $(x,x)$ of the diagonal. Since the diagonal is D-open in $\cl{R}_\varphi$, it follows that the upper triangle in the diagram above also commutes. As a consequence, $\ti{q}_* v = \ti{q'}_* v'$. Working across the entire chain, we conclude that $v_b = v_e$.
\end{proof}

\begin{proposition}
  \label{prop:groupoids-1}
If $\varphi\colon X \to Y$ is a Q-chart, then so is $T\varphi\colon TX \to TY$.
\end{proposition}
\begin{proof}
  Let $\varphi \colon X \to Y$ be a Q-chart. We use the following auxiliary result: the natural map $\nu$ determined by the diagram
  \begin{equation*}
    \begin{tikzcd}
      T(X \fiber{}{\varphi} X) \ar[drr, bend left, "T\beta"] \ar[dr, "\nu"] \ar[ddr, bend right, "T\alpha"] & & \\
      & TX \fiber{}{T\varphi} TX \ar[r] \ar[d] & TX \ar[d, "T\varphi"] \\
      & TX \ar[r, "T\varphi"] & FY
    \end{tikzcd}
  \end{equation*}
  is an isomorphism. First, observe that $TX$ and $TY$ are bundles over $X$ and $Y$, respectively, and $T(X \fiber{}{\varphi}X)$ and $TX \fiber{}{\varphi} TX$ are bundles over $X \fiber{}{\varphi} X$. Therefore, working fiberwise, we obtain the following commutative diagram, which we decorate using Lemma \ref{lem:groupoids-2} and the fact that $\alpha$ and $\beta$ are local diffeomorphisms.
    \begin{equation*}
    \begin{tikzcd}
      T_{(x,y)}(X \fiber{}{\varphi} X) \ar[drr, bend left, "T\beta", "\textrm{bij.}"'] \ar[dr, "\nu"] \ar[ddr, bend right, "T\alpha", "\textrm{bij.}"'] & & \\
      & T_xX \fiber{}{T\varphi} T_yX \ar[r] \ar[d] & T_xX \ar[d, "T\varphi", "\textrm{bij.}"'] \\
      & T_xX \ar[r, "T\varphi", "\textrm{bij.}"'] & T_yY
    \end{tikzcd}
  \end{equation*}
  The arrows labelled ``bij.''\ are bijections, and since these are bijections, every arrow in this diagram is a bijection. In particular, $\nu$ is a fiberwise bijection, and since it is a bundle map that is the identity on the base, $\nu$ is a bijection.

  The map $\nu$ is also an induction: for any D-open $A \subseteq X \fiber{}{\varphi} X$ on which both $\alpha$ and $\beta$ restrict to embeddings, we see that $\nu|_A \colon T(X \fiber{}{\varphi} X)|_A \to TX \fiber{}{T\varphi} TX$ is an induction. Finally, we note that any bijection that is locally an induction is a diffeomorphism.

  Having established that $\nu$ is a diffeomorphism, the fact that $T\varphi$ is a Q-chart readily follows:
  \begin{itemize}
  \item The map $T\varphi$ is a subduction because $T$ is a subduction.
  \item The diagonal $TX \fiber{}{1} TX$ is D-open in $TX \fiber{}{T\varphi} TX$ because we have
    \begin{equation*}
      \begin{tikzcd}
        T(X \fiber{}{1} X) \ar[r, "\subseteq", "\textrm{open}"'] \ar[d, "\cong"] & T(X \fiber{}{\varphi} X) \ar[d, "\nu", "\cong"'] \\
        TX \fiber{}{1} TX \ar[r, "\subseteq"] & TX \fiber{}{T\varphi} TX.
      \end{tikzcd}
    \end{equation*}
  \item The source and target for $TX \fiber{}{T\varphi} TX$ are local diffeomorphisms because they identify with $T\alpha$ and $T\beta$ under the diffeomorphism $\nu$.
  \end{itemize}
\end{proof}
The next two corollaries will be our main tools in the coming theorem.
\begin{corollary}
  \label{cor:groupoids-3}
  If $\varphi\colon X \to Y$ is a Q-chart, then so is $T^k\varphi$ for any $k \geq 0$. Furthermore
  \begin{itemize}
  \item whenever $T^k\varphi(v') = T^k\varphi(v)$, there is some $h \in \Psi(\varphi)$ with $T^kh(v') = v$;
  \item for all $x \in X$, the restriction $T_x^k\varphi\colon T^k_xX \to T^k_{\varphi(x)}Y$ is a bijection.
  \end{itemize} 
\end{corollary}
\begin{proof}
  That all the $T^k\varphi$ are Q-charts is a direct consequence of Proposition \ref{prop:groupoids-1}.

For the first item, we begin with the case $k = 1$. Suppose that $T\varphi(v) = T\varphi(v')$. Denoting the footprints of $v$ and $v'$ by $x$ and $x'$, we have $(x',x) \in \cl{R}_\varphi$, and so by Lemma \ref{lem:groupoids-3} we may find some $h \in \Psi(\varphi)$ with $h(x') = x$. Then $T\varphi(Th(v')) = T\varphi(v)$, and both $Th(v')$ and $v$ have footprint $x$, so by Lemma \ref{lem:groupoids-2} we have $Th(v') = v$. The general case follows inductively, invoking Proposition \ref{prop:groupoids-1}.

For the second item, surjectivity can be proved directly. Given $w \in T_y^kY$, we may find some $v' \in T^kX$ with $T\varphi(v') = w$, because $T\varphi$ is surjective. Denoting the footprint of $v'$ by $x'$, we may find $h \in \Psi(\varphi)$ with $h(x') = x$. Then $Th(v')$ is the element of $T_x^kX$ which $T\varphi$ sends to $w$.

For injectivity we proceed inductively. The case $k = 1$ is Lemma \ref{lem:groupoids-2}. Assuming the case for $k$, we show that $T_x^{k+1}\varphi$ is injective. Suppose that $T_x^{k+1}\varphi(v) = T_x^{k+1}\varphi(v')$. Applying $\pi_{T^kY}$ to both sides yields $T^k_x\varphi(\pi_{T^kX}(v)) = T^k_x\varphi(\pi_{T^kX}(v'))$, and by the inductive hypothesis $\pi_{T^kX}(v) = \pi_{T^kX}(v')$. But $T(T^k\varphi)$ is a Q-chart by Proposition \ref{prop:groupoids-1}, so by Lemma \ref{lem:groupoids-2} the equality $T_x^{k+1}\varphi(v) = T_x^{k+1}\varphi(v')$ implies that $v = v'$. 
\end{proof}

For $T_k\varphi$, we require fewer properties.
\begin{corollary}
  \label{cor:groupoids-2}
  If $\varphi \colon X \to Y$ is a Q-chart, then so is $T_k\varphi$ for any $k \geq 0$.
\end{corollary}
\begin{proof}
  First, we verify that $T_k\varphi$ is a subduction, so let $(p_i)\colon U \to T_kY$ be a plot, and fix $r \in U$. Denote $\pi(p_i(r)) = y$, which does not depend on $i$, and choose $x \in X$ with $\varphi(x) = y$. Since $T\varphi$ is a Q-chart, we may lift each $p_i$ near $r$ along $T\varphi$ to a plot $\ti{p_i}\colon (U,r) \to TX$, and moreover we can require that $\ti{p_i}(r) \in T_xX$.
  For each $i$, we have $\varphi \circ \pi \circ \ti{p}_i = \pi \circ p_1$, so all the $\pi \circ \ti{p}_i$ are lifts along $\varphi$ of the same map $\pi \circ p_i$. Since they all send $r$ to $x$, it follows that all the germs of $\pi \circ \ti{p_i}$ agree at $r$. Thus $(\ti{p_i})$ takes values in $T_kX$, at least near $r$, and its restriction the desired lift of $(p_i)$.

  Second, we verify that $T_kX \fiber{}{1} T_kX$ is D-open inside $T_kX \fiber{}{T_k\varphi} T_kX$. This follows from the fact that $TX \fiber{}{1} TX$ is D-open in $TX \fiber{}{T\varphi} TX$ (because $T\varphi$ is a Q-chart), and the fact we have the following two isomorphisms:
  \begin{align*}
    T_kX \fiber{}{1} T_kX &\cong (TX \fiber{}{1} TX) \fiber{}{\pi \times \pi} \cdots \fiber{}{\pi \times \pi} (TX \fiber{}{1} TX) \\
    T_kX \fiber{}{T_k\varphi} T_kX &\cong (TX \fiber{}{T\varphi} TX) \fiber{}{\pi \times \pi} \cdots \fiber{}{\pi \times \pi} (TX \fiber{}{T\varphi} TX)
  \end{align*}

  Finally, the second isomorphism above, and the fact that second projection $TX \fiber{}{T\varphi} TX \to TX$ is a local diffeomorphism (because $T\varphi$ is a Q-chart), together imply that the second projection $T_kX \fiber{}{T_k\varphi} T_kX \to T_kX$ is a local diffeomorphism.
\end{proof}

We are now equipped to treat the main result of this section, and the motivation for introducing Q-charts.
\begin{theorem}
  \label{thm:groupoids-2}
  Let $\varphi\colon X \to Y$ be a Q-chart. If $X$ is elastic, then so is $Y$.
\end{theorem}
\begin{proof}
For (E1), we begin with the diagram
    \begin{equation*}
      \begin{tikzcd}
        (\bb{L} \hat{T}_k)X \ar[r, "\theta_{k,X}", "\cong"'] \ar[d, "(\bb{L}\hat{T}_k)\varphi"] & T_kX \ar[d, "T_k\varphi"] \\
        (\bb{L}\hat{T}_k)Y \ar[r, "\theta_{k,Y}"] & T_kY.
      \end{tikzcd}
    \end{equation*}
    We want to show that $\theta_{k,Y}$ is an isomorphism. The top arrow is an isomorphism because $X$ satisfies (E1). The downward arrows are both subductions, the left one because $\varphi$ is a subduction, and the right one by Corollary \ref{cor:groupoids-2}. Therefore it suffices to show that $\theta_{k,Y}$ is injective.

    Fiberwise, both $\theta_{k,X}$ and $T_k\varphi$ are injective, the latter by Lemma \ref{lem:groupoids-2}. Therefore $\theta_{k,Y}$ is a fiberwise injection. But it is also a bundle map over the identity, therefore $\theta_{k,Y}$ is globally injective.

 For (E2), it suffices to find $\tau_Y$ making the diagram below commute.
    \begin{equation*}
      \begin{tikzcd}
        T^2X & T^2X \\
        T^2Y & T^2Y.
        \ar["\tau_X", from=1-1, to=1-2]
        \ar["T^2\varphi"', from=1-1, to=2-1]
        \ar["T^2\varphi", from=1-2, to=2-2]
        \ar["\tau_Y", from=2-1, to=2-2]
      \end{tikzcd}
    \end{equation*}
    The downward arrows are subductions, so all we need to show is that $\tau_X$ preserves the fibers of $T^2\varphi$. Suppose that $\nu,\nu' \in T^2X$ and that $T^2\varphi(\nu) = T^2\varphi(\nu')$. By Corollary \ref{cor:groupoids-3}, there is some $h \in \Psi(\varphi)$ with $T^2h(\nu') = \nu$. Then we invoke naturality of the canonical flip, namely that $\tau_X \circ T^2h = T^2h \circ \tau_X$:
    \begin{equation*}
      T^2\varphi \circ \tau_X (\nu') = T^2\varphi \circ \tau_X \circ T^2h (\nu) = T^2\varphi \circ T^2h \circ \tau_X(\nu) = T^2\varphi \circ \tau_X(\nu).
    \end{equation*}

    For (E3), we begin with the diagram
    \begin{equation*}
      \begin{tikzcd}
        TX & T^2X \\
        TY & T^2Y,
        \ar["\lambda_X", hook, from=1-1, to=1-2]
        \ar["T\varphi"', from=1-1, to=2-1]
        \ar["T^2\varphi", from=1-2, to=2-2]
        \ar["\lambda_Y", from=2-1, to=2-2]
      \end{tikzcd}
      \quad \textrm{and pass to the fibers} \quad
            \begin{tikzcd}
        T_xX & T^2_xX \\
        T_yY & T^2_yY.
        \ar["\lambda_X", hook, from=1-1, to=1-2]
        \ar["T_x\varphi"', "\cong", from=1-1, to=2-1]
        \ar["T^2_x\varphi", "\cong"', from=1-2, to=2-2]
        \ar["\lambda_Y", from=2-1, to=2-2]
      \end{tikzcd}
    \end{equation*}
    The maps $T_x\varphi$ and $T_x^2\varphi$ are bijections by Corollary \ref{cor:groupoids-2}. It immediately follows that $\lambda_Y$ is injective on the fibers, and since it is a bundle map that is identity on the base, we find that $\lambda_Y$ is globally injective.

    By a general diffeological argument, because $\lambda_X$ is an induction and $T\varphi$ and $T^2\varphi$ are subductions, we can conclude that the injective map $\lambda_Y$ is an induction if we show the containment
    \begin{equation*}
      (T^2\varphi)^{-1}(\lambda_Y(TY)) \subseteq \lambda_X(TX) .
    \end{equation*}
    We can show this fiber-wise: if $T^2\varphi(\nu) = \lambda_Y(w)$, then denoting the base-point of $\nu$ by $x$ and of $w$ by $y$, the diagram of fibers above lets us conclude that $\nu = \lambda_X((T_x\varphi)^{-1}(w))$.

    For (E4), to avoid cumbersome notation, we set
    \begin{equation*}
      (T^2X)_k \coloneqq T^2X \times_{TX}^{T\pi_X, T\pi_X} \cdots \times_{TX}^{T\pi_X, T\pi_X} T^2X.
    \end{equation*}
    Then we have the diagram
    \begin{equation*}
      \begin{tikzcd}
        TT_kX & (T^2X)_k & T_kTX \\
        TT_kY & (T^2Y)_k & T_kTY
        \ar["\nu_{k,X}", hook, from=1-1, to=1-2]
        \ar["{(\tau_X)_k}", "\cong"', from=1-2, to=1-3]
        \ar["TT_k\varphi", from=1-1, to=2-1]
        \ar["\nu_{k,Y}", from=2-1, to=2-2]
        \ar["{(\tau_Y)_k}", "\cong"', from=2-2, to=2-3]
        \ar["T_kT\varphi", from=1-3, to=2-3]
      \end{tikzcd}
    \end{equation*}
    The map $\tau_Y$ exists by (E2), which we have already showed. We can conclude that $\nu_{k,Y}$ is injective by working fiberwise, since $T_k\varphi$ and $T\varphi$ are fiberwise bijections by Lemma \ref{lem:groupoids-2} and Corollary \ref{cor:groupoids-2}.

   The axiom (E5) is a direct consequence of Corollaries \ref{cor:groupoids-3} and \ref{cor:groupoids-2}.
\end{proof}

\subsection{Elasticity of quotient groups}
\label{sec:elast-quot-groups}

Let $H$ be a normal subgroup of a diffeological group $G$ with inclusion $\iota\colon H \hookrightarrow G$. Then $G/H$ is a diffeological group, $T\iota(TH) \hookrightarrow TG$ is a normal subgroup, and by Corollary \ref{cor:1} we have an isomorphism of diffeological groups
\begin{equation*}
  T(G/H) \cong TG/T\iota(TH).
\end{equation*}
Whenever $T^k\iota \colon T^kH \to T^kG$ is an induction, we identify $T^kH$ with $T^k\iota(T^kH)$, and write $T^kG/T^kH$ for $T^kG/T^k\iota(T^kH)$.

\begin{maintheorem}{I}
\label{thm:groupoids-1}
Let $G$ be an elastic group, and let $H$ be a normal subgroup of $G$. If both $H$ and $G/H$ are elastic, and if iterated tangents $T^k\iota\colon T^kH \to T^kG$ of the inclusion $\iota\colon H \to G$ are inductions for all $k \geq 0$, then the natural map $\Lie(G)/\Lie(H) \to \Lie(G/H)$ is an isomorphism of Lie algebras.
\end{maintheorem}
\begin{proof}
Take the setup in the theorem's statement. Denote $\Lie(G)$ and $\Lie(H)$ by $\fk{g}$ and $\fk{h}$, respectively. We first need to verify that $\fk{h}$ is an ideal in $\fk{g}$. By Lemma \ref{lem:tangent-structure-1}, the Lie bracket on $\fk{g}$ is determined by $T_{(1)}T_{(2)}c$, where $c\colon G \times G \to G$ is the conjugation $(g,h) \mapsto ghg^{-1}$. By definition of a normal subgroup, whenever $h \in H$, we have $c(g,h) \in H$. Therefore $T_{(1)}T_{(2)}c$ restricts to a map
  \begin{equation*}
    T_{(1)}T_{(2)}c \colon TG \times TH \to T^2H.
  \end{equation*}
  Here we are using the assumption that $T^k\iota$ is an induction. This means that whenever $\xi \in \fk{g}$ and $\eta \in \fk{h}$, we have, by Lemma \ref{lem:tangent-structure-1}, that $\lambda_2(\eta, [\xi, \eta]) \in T^2H$. By definition $\lambda_2$, it follows that $[\xi, \eta] \in \fk{h}$, and so $\fk{h}$ is an ideal in $\fk{g}$.

  Now we show that $\eta_{G \rtimes H}\colon \Lie(G/H) \to \fk{g}/\fk{h}$ is an isomorphism of Lie algebras. Corollary \ref{cor:1} shows that $\eta_{G \rtimes H}$ is an isomorphism of diffeological spaces. It fits into the diagram
  \begin{equation*}
    \begin{tikzcd}
    & \fk{g} & \\
    \Lie(G/H) & & \fk{g}/\fk{h},
    \ar["T\pr"', two heads, from=1-2, to=2-1]
    \ar["\pr'", two heads, from=1-2, to=2-3]
    \ar["\eta_{G\rtimes H}", from=2-1, to=2-3]
  \end{tikzcd}
  \end{equation*}
  and the downward arrows are epimorphisms of Lie algebras: $T\pr$ because $\pr\colon G \to G/H$ is a surjective homomorphism between elastic groups, and $\pr'$ because we verified that $\fk{h}$ is an ideal in $\fk{g}$. Therefore $\eta_{G\rtimes H}$ is morphism of Lie algebras, and since it is bijective, it is an isomorphism of Lie algebras.
\end{proof}
Our first corollary is immediate.
\begin{corollary}
  \label{cor:groupoids-1}
  If $H$ is a diffeologically discrete normal subgroup of an elastic group $G$, then $H$ and $G/H$ are elastic and $\Lie(G/H) \cong \fk{g}$.
\end{corollary}
\begin{proof}
  Diffeologically discrete spaces are always elastic. The quotient $G/H$ is elastic because, by Proposition \ref{prop:groupoids-3}, the quotient map $G \to G/H$ is a Q-chart, and $G$ is assumed elastic. The maps $T^k\iota$ are inductions because $TH \cong H$. The fact $\Lie(G/H) \cong \fk{g}$ is then immediate.
\end{proof}

The next corollary is of independent interest, but we do not use it in this paper.
\begin{corollary}
  If $H$ is a normal subgroup of a finite-dimensional Lie group $G$, then $H$ and $G/H$ are elastic and $\Lie(G/H) \cong \fk{g}/\fk{h}$.
\end{corollary}
\begin{proof}
  By a theorem from Bourbaki \cite[III.4.5]{Bour89}, the subspace diffeology on $H$ is a manifold diffeology (i.e.\ $H$ is a Lie group), and so $H$ is elastic. The inclusion $\iota\colon H \to G$ is an induction, and $T\iota\colon TH \to TG$ is injective. In the trivializations, $T\iota$ has the form
  \begin{equation*}
    \begin{tikzcd}
      H \times \fk{h} & G \times \fk{g}
      \ar["{(\iota, T_0\iota)}", from=1-1, to=1-2]
    \end{tikzcd}
    \end{equation*}
    The map $\iota$ is an induction, and $T_0\iota\colon \fk{h} \to \fk{g}$ is a linear injection of finite-dimensional vector spaces equipped with their canonical diffeology, hence is also an induction. We conclude that $T\iota\colon TH \to TG$ is an induction, and inductively that $T^k\iota$ is an induction.

The fact that $G/H$ is elastic is a consequence of a result of Meigniez \cite[Corollary 5]{Meig97}, which states that the leaf space of any Riemannian foliation whose leaves are without holonomy is a Q-manifold. One can build the Q-charts directly by starting with an immersed submanifold $S \hookrightarrow G$ that meets every orbit of $H$ and such that $T_sG = T_sS \oplus T(s \cdot H)$ for each $s \in S$, and taking the quotient $S \to S/(s {\sim} s' \textrm{ if } s \in s' \cdot H) \cong G/H$. Therefore by Theorem \ref{thm:groupoids-2}, $G/H$ is elastic.
\end{proof}

\section{The convenient setting and diffeology}
\label{sec:conv-sett-diff}

Up to this point, our considerations and results have been fairly abstract, the only familiar example being quotients of finite-dimensional Lie groups. We now show that our setting includes the convenient analysis of Kriegl and Michor \cite{KriegMic97}: the category of convenient manifolds embeds into $\cat{Elst}$, and their tangent structures coincide. This allow us freely use all the constructions from the convenient calculus, and their examples become ours.

In Subsection \ref{sec:conv-diff-vect} we re-cast tangent-stability as convenience in the sense of Fr\"{o}licher and Kriegl \cite{FroelKrieg88}, and show that convenient spaces are elastic. In Subsection \ref{sec:from-conv-manif}, we embed Kriegl and Michor's category of convenient manifolds into $\cat{Elst}$. Doing so, we \emph{a fortiori} embed the categories of Banach and Fr\'{e}chet manifolds into $\cat{Elst}$. Crucially, we also show that their tangent structures coincide.

\subsection{Convenient diffeological vector spaces}
\label{sec:conv-diff-vect}

The diffeologies we encounter in this section will arise from collections of functions. Let $V$ be a set, and $\sr{F}_0$ be some collection of functions $V \to \R$. The set
\begin{equation*}
  \Pi \sr{F}_0 \coloneqq \{p\colon U \to V \mid U \in \cat{Cart}, \ f \circ p \textrm{ is smooth for all } f \in \sr{F}_0\}
\end{equation*}
is the diffeology on $V$ \define{determined by} $\sr{F}_0$.
\begin{definition}
  A diffeology $\sr{D}$ on $V$ is \define{reflexive} if $\sr{D} = \Pi C^\infty(V,\R)$.
\end{definition}
  The terminology ``reflexive,'' and the notation $\Pi$, is from \cite{BatKarWat23}.
  \begin{example}
    \label{ex:convenient-1}
  We have several examples.
  \begin{itemize}
  \item The diffeology $\Pi \sr{F}_0$ is always reflexive, and $\sr{F}_0 \subseteq C^\infty((V, \Pi \sr{F}_0),\R)$, but equality need not hold. For instance, the inclusion is strict if $\sr{F}_0$ contains only a single constant function.
  \item If $V = M$ is a finite-dimensional manifold, its canonical diffeology is reflexive.
  \item Consider the \define{spaghetti} diffeology on $\R^2$, which is generated by the smooth curves $C^\infty(\R, \R^2)$, i.e.\ it is finest diffeology admitting all smooth curves as plots. Observe that $C^\infty((\R^2, \sr{D}_{\text{standard}}), \R)$ and $C^\infty((\R^2, \sr{D}_{\text{spaghetti}}), \R)$ coincide. This is by Boman's theorem \cite{Bom67}, which grants that a function $\R^2 \to \R$ is smooth if and only if its pullback by all smooth curves is smooth. Thus both collections of functions determine the standard diffeology on $\R^2$, and the spaghetti diffeology must not be reflexive.
  \end{itemize}
\end{example}

Reflexive diffeological spaces form a full subcategory of $\cat{Dflg}$.

\begin{remark}
  Given a reflexive diffeological space $(V, \sr{D})$, we can form the triple
  \begin{equation*}
  (V, C^\infty(\R,V), C^\infty(V,\R)).  
\end{equation*}
The data of this triple is an example of a \define{Fr\"{o}licher space}. In fact, by \cite[Theorem 2.13]{BatKarWat23} this assignment of a Fr\"{o}licher space to each reflexive diffeological space is an isomorphism between the categories of reflexive diffeological spaces, and Fr\"{o}licher spaces.  Fr\"{o}licher and Kriegl \cite{FroelKrieg88} work with the category of Fr\"{o}licher spaces, which they call $\cl{L}\textit{ip}^\infty$-spaces; the ``Lip'' stands for ``Lipschitz.''
\end{remark}

\begin{notation}
  We will denote by $\cat{DflgVS}$ the category of diffeological vector spaces, and smooth linear maps between them. We will write ``dvs'' for ``diffeological vector space.'' We write $L^\infty(V,W) \coloneqq \underline{\cat{DflgVS}}(V,W)$, and if $W = \R$, then we simply write $L^\infty(V)$.
\end{notation}
Given a diffeological vector space $V$, and a smooth curve $c\colon \R \to V$, we would like to define its velocity curve $\dot{c}\colon \R \to V$. For $V = \R^n$, recall that the tangent of $c$ is defined by
\begin{equation*}
\hat{T}c(t,t_0) = (c(t), \dot{c}(t)t_0),
\end{equation*}
where $\dot{c}(t)$ is the usual velocity vector. Thus $\dot{c}(t)$ is the second component of $ \hat{T}c(t,1)$. For a general dvs $V$, we no longer necessarily have $V \times V \cong TV$. But we have the natural map $\ell$ from Section \ref{sec:elast-diff-groups-1}:
\begin{equation*}
\ell\colon V \times V \to TV, \quad (u,v) \mapsto \partial[u+tv]. 
\end{equation*}
Recall that we call $V$ tangent-stable if $\ell$ is an isomorphism.
\begin{definition}
  Given a smooth curve $c\colon \R \to V$ into a tangent-stable dvs, its \define{velocity curve}, denoted $\dot{c}\colon \R \to V$, is the smooth curve defined by the commutative diagram
  \begin{equation*}
    \begin{tikzcd}
      \R \ar[r, "\dot{c}"] \ar[d, "{(\cdot, 1)}"] & V \\
      \R \times \R \ar[d, "\ell", "\cong"'] \ar[r] & V \times V \ar[u, "\pr_2"] \ar[d, "\ell", "\cong"'] \\
      T\R \ar[r, "Tc"] & TV.
    \end{tikzcd}
  \end{equation*}
\end{definition}

\begin{proposition}
\label{prop:4}
Fix a smooth curve $c\colon \R \to V$ into a tangent-stable dvs. For each $l \in L^\infty(V)$, the composition $l \circ c$ is smooth, and $\dot{\wideparen{l \circ c}} = l \circ \dot{c}$.
\end{proposition}
\begin{proof}
  Smoothness of $l \circ c$ is clear. Unravelling definitions, and using the chain rule, we write
  \begin{align*}
    \dot{\wideparen{l \circ c}}(t) &= \pr_2 \circ \ell^{-1} \circ T(l \circ c) \circ \ell (t,1) \\
    &= \pr_2 \circ (\ell^{-1} \circ Tl \circ \ell) \circ (\ell^{-1} \circ Tc \circ \ell)(t,1).
  \end{align*}
  But because $l$ is linear, $Tl \circ \ell = \ell \circ (l \times l)$. For clarity, we confirm this:
  \begin{equation*}
    Tl(\ell(u,v)) = Tl(\partial[u+tv]) = \partial [l(u) + tl(v)] = \ell(l(u),l(v)).
  \end{equation*}
  Thus
  \begin{equation*}
    \pr_2 \circ (\ell^{-1}\circ  Tl \circ \ell) = \pr_2 \circ (l \times l) = l \circ \pr_2.
  \end{equation*}
  We conclude
  \begin{equation*}
    \dot{\wideparen{l \circ c}}(t) = l \circ (\pr_2 \circ \ell^{-1} \circ Tc \circ \ell)(t_0,1) = l \circ \dot{c}(t).
  \end{equation*}
\end{proof}

Proposition \ref{prop:4} says that the velocity curve of a smooth curve is also its ``weak derivative,'' which is a notion we can define for a general dvs.

\begin{definition}
A curve $c\colon \R \to V$ into a dvs is \define{weakly differentiable} (with respect to $L^\infty(V)$), with a \define{weak derivative} $\dot{c}\colon \R \to V$, if for all $l \in L^\infty(V)$, the composition $l \circ c$ is differentiable, and its derivative satisfies $\dot{\wideparen{l\circ c}} = l \circ \dot{c}$.
\end{definition}

A dvs may, or may not, satisfy the following properties:

\begin{itemize}
\item infinitely weakly differentiable curves are smooth;
\item weak derivatives, when they exist, are unique;
\item smooth curves are infinitely weakly differentiable.
\end{itemize}
Each property is a consequence of a natural condition imposed on $V$. We describe these conditions sequentially, and we will call those dvs that satisfy all three ``convenient.'' We have proved that tangent-stable dvs satisfy the third property.

\begin{definition}
  \label{def:convenient-1}
  We introduce three terms:
  \begin{itemize}
  \item   A diffeological vector space $(V, \sr{D}_V)$ is \define{pre-convenient} if $\sr{D}_V = \Pi L^\infty(V)$.
  \item   A diffeological vector space $V$ is \define{linearly separated} if $L^\infty(V)$ separates points.
  \item   A pre-convenient and linearly separated dvs $V$ is called \define{convenient} if every smooth curve is infinitely weakly differentiable.
  \end{itemize}
\end{definition}

Pre-convenient dvs are in particular reflexive, but we are not sure if the converse holds. The terminology ``pre-convenient'' comes from \cite{FroelKrieg88}.

\begin{lemma}
  Two items:
  \begin{itemize}
  \item  In a pre-convenient dvs, infinitely weakly differentiable curves are smooth.
  \item In a linearly separated dvs, weak derivatives are unique.
  \end{itemize}    
\end{lemma}
\begin{proof}
  Item by item:
  \begin{itemize}
  \item  This is clear from the definition: since $\sr{D}_V = \Pi L^\infty(V)$, smoothness of $c\colon \R \to V$ is equivalent to smoothness of $l \circ c$ for all $l \in L^\infty(V)$, which holds for infinitely weakly differentiable curves.
  \item   This is also clear from the definition: if $\dot{c}$ and $d$ are two weak derivatives of $c$, then for all $l \in L^\infty(V)$,
  \begin{equation*}
    l \circ \dot{c}(t_0) = \dot{\wideparen{l \circ c}}(t_0) = l \circ d(t_0),
  \end{equation*}
  thus by linear separation, $\dot{c}(t_0) = d(t_0)$.
  \end{itemize}
 \end{proof}
  \begin{remark}
    This definition of a convenient dvs $V$ is precisely the one used in \cite{FroelKrieg88}, where they view $V$ as a Fr\"{o}licher space. Specifically, in their notation, it is condition (5) in \cite[Theorem 2.6.2]{FroelKrieg88}, with $k=\infty$.
  \end{remark}

  We proved in Proposition \ref{prop:4} that pre-convenient and linearly separated dvs that are tangent-stable are convenient. Now we treat the converse. We will need the notion of a weak integral.

\begin{definition}
  A curve $c\colon \R \to V$ into a dvs is \define{weakly integrable}, with a \define{weak integral} $\int c\colon \R \to V$, if for all $l \in L^\infty(V)$, the composition $l \circ c$ is integrable, and its integral satisfies $\int_0^t l \circ c = l \circ \int c(t)$.
\end{definition}

\begin{proposition}[{\cite[Theorem 2.6.2]{FroelKrieg88}}]
  If $V$ is a convenient dvs, then every smooth curve is weakly integrable, and its weak integral is unique.
\end{proposition}

  This is a fundamental result in the theory of convenient calculus. In fact, on a pre-convenient and linearly separated dvs, weak differentiability of smooth curves is equivalent to weak integrability.

\begin{lemma}[Fundamental theorem of calculus]
  \label{lem:4}
  Fix a convenient dvs $V$. The maps
  \begin{equation*}
    C^\infty(\R, V) \to C^\infty(\R, V), \quad c \mapsto \dot{c} \text{ and } c \mapsto \int c
  \end{equation*}
  are smooth and linear. Furthermore, $t \mapsto c(t)-c(0)$ is the weak integral of $\dot{c}$. 
\end{lemma}
We will denote $\int_s^tc \coloneqq \int c(t) - \int c(s)$. Thus $\int c(t) = \int_0^t c$.
\begin{proof}
  Take any plot $p\colon U \to C^\infty(\R,V)$. We handle differentiation first, so we must show that $(r,t) \mapsto \dot{p_r}(t)$ is smooth. By pre-convenience of $V$, we may check this in each $l \in L^\infty(V)$: we must show that $(r,t) \mapsto l \circ \dot{p_r}(t) = \dot{\wideparen{l \circ p_r}}(t)$ is smooth. But this map is simply the partial derivative of the smooth map $(r,t) \mapsto l \circ p_r(t)$ in the $t$-direction, hence is itself smooth.

  For integration, we must show that $(r,t) \mapsto \int p_r(t)$ is smooth. By pre-convenience of $V$, we may check this in each $l \in L^\infty(V)$: we must show that $(r,t) \mapsto l \circ \int p_r(t) = \int_0^t l \circ p_r$ is smooth. But this map is simply the integration of the smooth map $(r,t) \mapsto l \circ p_r(t)$, hence is itself smooth.

  For the remaining assertions, check all required equalities scalar-wise.
\end{proof}

Given a smooth map $S\colon \R^n \to V$ into a convenient dvs, we have its partial derivatives $\partial_iS \colon \R^n \to V$, defined by
\begin{equation*}
  \begin{tikzcd}
    & C^\infty(\R,V) & \\
    \R^n & & V. \\
    \ar["{\vec{t} \ \mapsto\  S(\vec{t} + se_i)}", from=2-1, to=1-2]
    \ar["{c \ \mapsto\  \dot{c}(0)}", from=1-2, to=2-3]
    \ar["\partial_iS", from=2-1, to=2-3]
  \end{tikzcd}
\end{equation*}

\begin{lemma}[Hadamard's lemma]
  \label{lem:5}
  Fix a convenient dvs $V$. If $S\colon \R^n \to V$ is smooth, then there exists smooth maps $h_i\colon \R^n \to V$ with $c(\vec{t}) = c(0) + \sum t_ih_i(\vec{t})$. Necessarily $h_i(0) = \partial_iS(0)$.
\end{lemma}
\begin{proof}
  Consider the auxiliary smooth curve $H(s) \coloneqq c(s\vec{t})$. Then $\dot{H}(s) = \sum t_i \partial_iS(s\vec{t})$ (check scalar-wise, and use the chain rule for functions $\R^n \to \R$), and so by the fundamental theorem of calculus and linearity of the integral,
  \begin{equation*}
    c(t) - c(0) = H(1)-H(0) = \int_0^1 \sum_i t_i\partial_iS(s\vec{t}) ds = \sum t_i \int_0^1 \partial_iS(s\vec{t})ds.
  \end{equation*}
  Thus setting $h_i(t) \coloneqq \int_0^1\partial_iS(s\vec{t})ds$ will work, if $h_i$ is smooth. This is true because $h_i$ is the composition of the smooth maps:
  \begin{equation*}
    \begin{tikzcd}
      \R^n & V \\
      C^\infty(\R,V) & C^\infty(\R,V)
      \ar["{\vec{t} \ \mapsto \ [S(s\vec{t})]}"', from=1-1, to=2-1]
      \ar["h", from=1-1, to=1-2]
      \ar["\int", from=2-1, to=2-2]
      \ar["c \ \mapsto \ c(1)"', from=2-2, to=1-2]
    \end{tikzcd}
  \end{equation*}
\end{proof}

\begin{proposition}
  \label{prop:5}
  Every convenient vector space is tangent-stable.
\end{proposition}

\begin{proof}
  We must show that $\ell\colon V \times V \to TV$ is an isomorphism of diffeological spaces. We will use the maps
  \begin{alignat*}{2}
    (\ev_0,D_0)&\colon C^\infty(\R,V) \to V\times V \quad & c &\mapsto (c(0),\dot{c}(0)) \\
  \gamma&\colon V \times V \to C^\infty(\R,V) \quad & (u,v) &\mapsto [u+tv] .
  \end{alignat*}
  These are both smooth, $D_0$ because it is the composition of the evaluation map $\ev_0$ with the smooth map $c \mapsto \dot{c}$, and $\gamma$ because the map $(u,v,t) \mapsto \gamma(u,v)(t)$ is smooth.

  We claim that $\gamma$ is a section of $(\ev_0,D_0)$. Clearly $\ev_0(\gamma(u,v)) = v$. As for $D_0(\gamma(u,v))$, observe that for each $l \in L^\infty(V)$,
  \begin{equation*}
        l \circ D_0(\gamma(u,v)) = \frac{d}{dt}\Big|_{t=0}\left(l(u) + tl(v)\right) = l(v),
      \end{equation*}
      and since $V$ is linearly separated, $D_0(\gamma(u,v)) = v$. Having a section, the map $(\ev_0,D_0)$ is a subduction.

      Now we will show the diagram below commutes:
      \begin{equation}
        \label{eq:convenient-1}
    \begin{tikzcd}
      & C^\infty(\R, V) & \\
      V \times V  & & TV
      \ar["{(\ev_0,D_0)}"', two heads, from=1-2, to=2-1]
      \ar["\partial", two heads, from=1-2, to=2-3]
      \ar["\ell", from=2-1, to=2-3]
    \end{tikzcd}
  \end{equation}
Fix a curve $c$, and denote $(u,v) \coloneqq (c(0), \dot{c}(0))$. Applying Hadamard's lemma \ref{lem:5} to $c$ twice, we find smooth $h$ with $c(t) = u + tv + t^2h(t)$. Then, setting $H(t,t') \coloneqq u + tv +t'h(t)$, we get the commutative diagram
    \begin{equation*}
      \begin{tikzcd}[sep=large]
        \R & \R^2 & \R \\
        & V &
        \ar["{t\ \mapsto \ (t,t^2)}", from=1-1, to=1-2]
        \ar["{(t,0)\ \mapsfrom \ t}"', from=1-3, to=1-2]
        \ar["c", from=1-1, to=2-2]
        \ar["H", from=1-2, to=2-2]
        \ar["{\gamma(v,v_0)}", from=1-3, to=2-2]
      \end{tikzcd}
    \end{equation*}
    This diagram shows that $\partial c = \partial [u+tv]$, as required.

Since Diagram \eqref{eq:convenient-1} commutes, and the downward arrows are subductions, $\ell$ is also a subduction. To conclude that $\ell$ is an isomorphism, it suffices to show that it is injective. Suppose that $\partial[u_b + tv_b] = \partial[u_e+tv_e]$ (the subscripts denote ``beginning'' and ``end''). An arbitrary link in the chain witnessing this equality has the form (in Notation \ref{not:tangent-structure-2})
    \begin{equation*}
      \begin{tikzcd}
        (U, v) \ar[rr, "f"] \ar[dr, "p"'] & & (U', v') \ar[dl, "p'"] \\
         & V. &
      \end{tikzcd}
    \end{equation*}
    Let $c$ and $c'$ denote any curves in $U$ and $U'$ whose derivatives at $0$ are $v$ and $v'$, respectively. Then, for any $l \in L^\infty(V)$, observe:
    \begin{equation*}
      l(\dot{\wideparen{pc}}(0)) = \dot{\wideparen{lpc}}(0) = \dot{\wideparen{lp'fc}}(0) = \hat{T}(lp') \circ \dot{\wideparen{fc}}(0) = \hat{T}(lp') \circ \dot{\wideparen{c'}}(0),
    \end{equation*}
    because $\dot{\wideparen{fc}}(0) = \dot{\wideparen{c'}}(0) = v'$. Thus $\dot{\wideparen{pc}}(0) = \dot{\wideparen{p'c'}}(0)$.

    At the beginning of the chain, clearly $p(u_b+tv_b)$ has weak derivative $v_b$ at $t=0$, and similarly $p(u_e + tv_e)$ has weak derivative $v_e$. Since these weak derivatives are preserved across the chain, we find that $v_b = v_e$. As $u_b = u_e$ a necessary condition for the equality, we have proved that $\nu$ is injective.
  \end{proof}

  We have now arrived at a statement that summarizes the results of this section.

  \begin{corollary}
    \label{cor:convenient-2}
    A pre-convenient and linearly separated dvs is convenient if and only if it is tangent-stable.
  \end{corollary}
  
   \begin{proof}
 This directly follows from Propositions \ref{prop:4} and \ref{prop:5}.
\end{proof}

Finally, we relate these notions to elasticity.

   \begin{proposition}
     \label{prop:convenient-2}
     A convenient dvs is elastic.
   \end{proposition}
   \begin{proof}
     Fix a convenient dvs $V$. By Proposition \ref{prop:5}, we know that $V$ is tangent-stable. To show that $V$ is elastic, we appeal to Lemma \ref{lem:tangent-structure-7}. See Subsection \ref{sec:elast-diff-groups-1} for the relevant notation.

    By Lemma \ref{lem:tangent-structure-7}, we only need to show that $N^2$ and the $N_k$ are surjective. For surjectivity of $N^2$, observe that an arbitrary element $p_*(u,v_0,v_1,v_{01})$ of $(\bb{L}\hat{T}^2)V$ can be represented as $S_*(0,e_1,e_2,0)$ for some parametrized surface $S$: to do this, pull $p$ back along a map $\R^2 \to U= \dom p$ whose germ at the origin coincides with $(s,t) \mapsto u + sv_0 + tv_1 + stv_{01}$. Apply Hadamard's lemma to $S$ several times, (re)-instantiating the vectors $u,v_0,v_1,v_{01},v_{10} \in V$, as follows:
     \begin{align*}
       S(s,t) &= u + sS_0 + tS_1 \\
                   &= u + s(v_0 + sS_{00} + tS_{01}) + t(v_1 + sS_{10} + tS_{11}) \\
                   &= u + sv_0 + tv_1 \\
                   &\quad + st(v_{01} + sS_{010} + tS_{011}) + ts(v_{10} + sS_{100} + tS_{101}) \\
                   &\quad + s^2S_{10} + t^2S_{11} \\
       &= u + sv_0 + tv_1 + st(v_{01}+v_{10}) + s^2h_0(s,t) + t^2h_1(s,t).
     \end{align*}
Setting $H(s,t,s',t') \coloneqq u + sv_0 + tv_1 + st(v_{01}+v_{10}) + s'h_0 + t'h_1$, and setting $\vec{v} \coloneqq (u,v_0,v_1,v_{01}+v_{10})$, we get
     \begin{equation*}
       \begin{tikzcd}[column sep=4.0cm]
         \R^2 & \R^2 \times \R^2 & \R^2 \\
         & V, &
         \ar["{(s,t) \ \mapsto \ (s,t,s^2,t^2)}", from=1-1, to=1-2]
         \ar["{(s,t,0,0)\ \mapsfrom\ (s,t)}"', from=1-3, to=1-2]
         \ar["{S}", from=1-1, to=2-2]
         \ar["{S^2(\vec{v})}", from=1-3, to=2-2]
         \ar["H", from=1-2, to=2-2]
       \end{tikzcd}
     \end{equation*}
     and this diagram witnesses the equality $S_*(0,e_1,e_2,0) = S^2(\vec{v})_*(0,e_1,e_2,0) = N^2(\vec{v})$.

     To show that $N_k$ is onto, observe that an arbitrary element $p_*(u,v_0,\ldots, v_k)$ of $(\bb{L}\hat{T}_k)V$ can be represented as $S_*(0,e_1,\ldots, e_k)$ for some $S \colon \R^k \to V$; to do this, pull $p$ back along a map $\R^k \to U = \dom p$ whose germ at the origin coincides with $\vec{t} \mapsto u + \sum t_i v_i$. Apply Hadamard's lemma to $S$ twice, (re)-instantiating the vectors $u,v_0,\ldots, v_k \in V$ as follows:
     \begin{align*}
       S(\vec{t}) &= u + \sum t_i S_i \\
                       &= u + \sum t_i\left(v_i + \sum_j t_j S_{ij}\right) \\
       &= u + \sum t_i v_i + \sum_{i,j} t_it_j S_{ij}.
     \end{align*}
     Setting $H(\vec{t},\vec{r}) \coloneqq u + \sum_i t_iv_i + \sum_{i,j} r_{ij}S_{ij}$, where $\vec{r} \in \R^{\frac{(k-1)k}{2}}$, and setting $\vec{v} \coloneqq (u,v_0\ldots, v_{k-1})$, we get
     \begin{equation*}
       \begin{tikzcd}[column sep=4.0cm]
         \R^k & \R^k \times \R^{\frac{(k-1)k}{2}} & \R^k \\
         & V &
         \ar["{\vec{t} \ \mapsto \ (\vec{t}, \ (t_it_j)_{ij})}", from=1-1, to=1-2]
         \ar["{(\vec{t},0)\ \mapsfrom\ \vec{t}}"', from=1-3, to=1-2]
         \ar["{S}"', from=1-1, to=2-2]
         \ar["{S_2(\vec{v})}", from=1-3, to=2-2]
         \ar["H", from=1-2, to=2-2]
       \end{tikzcd}
     \end{equation*}
     and this diagram witnesses the equality $S_*(0,e_1,\ldots, e_k) = S_2(\vec{v})_*(0,e_1,\ldots, e_k)$.
   \end{proof}

   \begin{remark}
Suppose that $V$ is pre-convenient and linearly separated. Then it is convenient if and only if it is tangent-stable, in which case it is furthermore elastic. We did not show that if $V$ is elastic, then it is tangent-stable.
   \end{remark}

   \begin{remark}
     \label{rem:convenient-1}
     When $V$ is tangent-stable and elastic, and $U \subseteq V$ is D-open, we have $TU \cong U \times V$, and the rest of the tangent structure is identical to that on $\cat{Cart}$. This also helps describe the tangent structure maps on any space that is locally diffeomorphic to $V$.
   \end{remark}

\subsection{From convenient manifolds to diffeology}
\label{sec:from-conv-manif}

Now we import topological vector spaces into diffeology. We recall that a \define{locally convex topological vector space}, abbreviated ``lctvs,'' is a vector space $E$ equipped with a vector space topology that is generated by a family of seminorms. A \define{morphism} of lctvs is a continuous linear map.

\begin{notation}
We denote by $\cat{LCTVS}$ the category of lctvs and continuous linear maps between them. We will write $L^*(E, E_1) \coloneqq \underline{\cat{LCTVS}}(E, E_1)$, and when $E_1 = \R$, we write simply $L^*(E)$.
\end{notation}

\begin{definition}
  A lctvs is \define{separated} if the seminorms defining its topology separate points.
\end{definition}

\begin{example}
  A \define{Banach} space is a lctvs whose topology is determined by a complete norm. A \define{Fr\'{e}chet} space is a separated lctvs whose topology is generated by a countable collection of seminorms, and which is complete with respect to these seminorms.
\end{example}

Given a separated lctvs $E$, we consider two reasonable diffeologies.
\begin{itemize}
\item The diffeology $\sr{D}_{\text{c}}(E)$. Given a map $c\colon U \to E$, with $U \subseteq \R$ open, we call $c$ \define{differentiable} if for each $t \in U$, the difference quotient
  \begin{equation*}
    c' (t) = c^{(1)}(t) \coloneqq \lim_{h \to 0} \frac{c(t+h)-c(t)}{h}
  \end{equation*}
  exists. We inductively define $c^{(k)}$, and declare $c$ to be $\sr{D}_{\text{c}}(E)$\define{-smooth}\footnote{the ``c'' stands for ``convenient''} if $c^{(k)}$ exists for each $k \geq 0$. The plots of $\sr{D}_{\text{c}}(E)$ are those parametrizations $p\colon U \to E$ for which $p \circ c$ is $\sr{D}_{\text{c}}(E)$-smooth for all smooth (in the usual sense) curves $c\colon V \to U$.

  \item The diffeology $\Pi L^*(E)$.
  \end{itemize}

  In the diffeology $\sr{D}_{\text{c}}(E)$, it is apparent and natural that the velocity of the curve $c$ is $c'$. However, the assignment $E \mapsto (E, \sr{D}_{\text{c}}(E))$ is not a functor from $\cat{LCTVS}$ to $\cat{Dflg}$, because continuous linear maps between lctvs need not be $\sr{D}_{\text{c}}$-smooth. On the other hand, the assignment $E \mapsto \Pi L^*(E)$ is functorial. We always have $\sr{D}_{\text{c}}(E) \subseteq \Pi L^*(E)$.

  \begin{definition}[Kriegl-Michor \cite{KriegMic97}]
    A separated lctvs $E$ is \define{convenient} if $\sr{D}_{\text{c}}(E) = \Pi L^*(E)$. 
  \end{definition}

  The condition $\sr{D}_{\text{c}}(E) = \Pi L^*(E)$ is not sufficient to ensure that the diffeological vector space $(E, \Pi L^*(E))$ is convenient in the sense of Definition \ref{def:convenient-1}. This is because, in general, $L^*(E) \subseteq L^\infty(E, \Pi L^*(E))$, and equality need not hold. A remedy is to consider only bornological lctvs.

  \begin{definition}
    A lctvs $E$ is \define{bornological} if every bounded seminorm on $E$ is continuous.
  \end{definition}

  \begin{remark}
    The seminorms that define the topology of $E$ are by definition bounded and continuous. However, an arbitrary seminorm on $E$ may not be bounded, and if it is bounded, it may not be continuous. In other words, the definition of bornological lctvs is not circular, and is not a tautology.
  \end{remark}
  
  \begin{example}
    Banach and Fr\'{e}chet spaces are bornological.
  \end{example}
  
  \begin{maintheorem}{II}
\label{thm:convenient-1}
  A separated bornological locally convex topological vector space is convenient if and only if, when viewed as a diffeological space, it is tangent-stable.
  \end{maintheorem}

  \begin{proof}
    Fix a separated and bornological lctvs $E$. We will show the equivalence of the following three properties:
    \begin{itemize}
\item $(E, \Pi L^*(E))$ is convenient.
\item $(E, \Pi L^*(E))$ is tangent-stable.
\item $\sr{D}_{\text{c}}(E) = \Pi L^*(E)$ 
\end{itemize}
    Suppose that $(E, \Pi L^*(E))$ is convenient. Then every smooth curve $c$ has a weak-derivative $\dot{c}$, and in fact $\dot{c} = c'$.  This shows that the 1-plots of $\Pi L^*(E)$ agree with the 1-plots of $\sr{D}_{\text{c}}(E)$. But both diffeologies are determined by their 1-plots, so the diffeologies agree.

    Conversely, if $\sr{D}_{\text{c}}(E) = \Pi L^*(E)$, then given any smooth curve $c$, we may always form $c'$. For every $l \in L^*(E)$, we have $\dot{\wideparen{l \circ c}} = l \circ c'$. Then, because $E$ is bornological, $L^\infty(E, \Pi L^*(E)) = L^*(E)$, so $c'$ is the required weak-derivative of $c$.

    Since $E$ is Hausdorff and bornological, convenience of $(E, \Pi L^*(E))$ is equivalent to tangent-stability of $(E, \Pi L^*(E))$ by Proposition \ref{cor:convenient-2}.
  \end{proof}
We now work towards convenient manifolds. We need the notion of a non-linear smooth map between convenient vector spaces. For two separated lctvs $E$ and $F$, Kriegl and Michor take a smooth map between them to be any $f\colon E \to F$ which is $\sr{D}_{\text{c}}$-smooth. Thus if $E$ and $F$ are convenient lctvs, the smooth maps between them are precisely the diffeologically smooth ones, when they are viewed as convenient diffeological vector spaces.
  
  \begin{remark}
    \label{rem:convenient-2}
    There is another notion of smooth map between lctvs. Let $E$ and $F$ be separated lctvs. We call a map $f\colon E \to F$ \emph{Michel-Bastiani smooth} if all its iterated directional derivatives
    \begin{equation*}
      Df\colon E \times E \to F, \quad D^2f\colon (E \times E) \times (E \times E) \to F, \ldots,
    \end{equation*}
    exist. Here $Df(u,v) \coloneqq [t \mapsto f(u+tv)]'(0)$. Every MB-smooth map is $\sr{D}_{\text{c}}$-smooth, but the converse need not hold, even if $E$ and $F$ are convenient lctvs; Gl\"{o}ckner \cite{Gloec06} gives a counter-example. However, when $E$ and $F$ are Fr\'{e}chet spaces, Fr\"{o}licher \cite{Froel81} proved that every $\sr{D}_{\text{c}}$-smooth map is Michel-Bastiani smooth. 
  \end{remark}
  
Now we need a notion of locality. Fix a convenient lctvs $E$. In defining manifolds modelled on $E$, Kriegl and Michor equip $E$ not with its original topology, but with its ``$\text{c}^\infty$-topology.'' This is the final topology with respect to the c-smooth curves. This is the same as the $D$-topology of $(E, \sr{D}_{\text{c}}(E))$, thanks to the following result.
\begin{lemma}[{\cite[Theorem 3.7]{ChrSinWu14}}]
  \label{lem:convenient-2}
  Let $X$ be a diffeological space. Its $D$-topology is the final topology with respect to its smooth curves $C^\infty(\R, X)$.
\end{lemma}

Kriegl and Michor proceed to define a convenient manifold $M$ modelled on $E$ precisely as one does for finite-dimensional manifolds modelled on $\R^n$: charts take values in $\text{c}^\infty$-open subsets of $E$, and transition functions are required to be $\sr{D}_{\text{c}}$-smooth. Smooth maps between convenient manifolds are those maps which are $\sr{D}_{\text{c}}$-smooth in charts. They equip $M$ with the final topology induced by the charts, and require this to be Hausdorff. This yields a category of manifolds $\cat{ConvMan}$.

\begin{example}
Banach and Fr\'{e}chet manifolds are convenient manifolds. This is because the original topology on a Banach or Fr\'{e}chet space, and their $\text{c}^\infty$-topology, agree, and because the Michel-Bastiani smooth maps coincide with the $\sr{D}_{\text{c}}$-smooth ones (cf.\ Remark \ref{rem:convenient-2}).
\end{example}

On the other hand, we may also consider those diffeological spaces which are locally diffeomorphic to convenient dvs, and whose D-topology is Hausdorff. Exactly as in the case of finite-dimensional manifolds, there is no difference in the two approaches.

\begin{proposition}
  \label{prop:convenient-3}
  Let $\cat{ConvMan}$ denote the category of convenient manifolds as described above. The canonical functor
  \begin{equation*}
    \cat{ConvMan} \to \cat{Dflg},
  \end{equation*}
  which assigns to a convenient manifold $M$ the diffeology consisting of $\sr{D}_{\text{c}}$-smooth maps $U \to M$, is full and faithful. Its essential image consists of those diffeological spaces which are locally diffeomorphic to convenient spaces, and which have Hausdorff D-topology.
\end{proposition}
\begin{proof}
  Proceed as in \cite{Igl13}, which treats the finite-dimensional case.
\end{proof}

\begin{remark}
  We have seen that convenient lctvs are all pre-convenient, hence reflexive, diffeological vector spaces. However, reflexivity is a property that any diffeological space may possess. It is not clear if convenient manifolds embed into the subcategory of reflexive diffeological spaces.
\end{remark}

This embedding is also compatible with the relevant tangent functors. Kriegl and Michor define the tangent functor $\hat{S}$ on $\text{c}^\infty$-open subsets of convenient lctvs by
\begin{equation*}
  \hat{S}(U \subseteq E) \coloneqq U \times E, \quad \hat{S}f(u,v) \coloneqq (f(u), [t \mapsto f(u + tv)]'(0)),
\end{equation*}
and then define $S\colon \cat{ConvMan} \to \cat{ConvMan}$ as one usually does for finite-dimensional manifolds.

\begin{proposition}
  \label{prop:convenient-4}
  When we view $\cat{ConvMan}$ as a subcategory of $\cat{Elst}$, its tangent structure coincides with ours on $\cat{Elst}$.
\end{proposition}
\begin{proof}
  Fix a convenient manifold $M$, and a some D-open subset $U$ which is diffeomorphic to a D-open subset of a convenient space $V$.  Remark \ref{rem:convenient-1} details the tangent structure on $U$. But all of these maps coincide with those from tangent structure used by \cite{KriegMic97}, because on open subsets of $V$, we have $\hat{S} = T$. Thus we have the same tangent structure.
\end{proof}

From this discussion, it follows that any results or computations for convenient, and in particular Banach or Fr\'{e}chet, manifolds from the literature are also valid in $\cat{Dflg}$. For instance, if $G$ is a Banach-Lie group, meaning a group object in the category of convenient manifolds, modelled on a Banach space, with Banach-Lie algebra $\fk{g}$, then when we view $G$ as a diffeological group, it is elastic with Lie algebra $\fk{g}$.

\section{Integrating central extensions}
\label{sec:integr-centr-extens}

In this final section, we integrate some non-enlargeable Banach-Lie algebras to elastic diffeological groups. In Subsection \ref{sec:diff-centr-extens} we review central extensions of diffeological groups and diffeological Lie algebras, culminating in a differentiation procedure for central extensions of elastic groups. In Subsection \ref{sec:integr-conv-lie}, we restrict our attention to central extensions of convenient Lie algebras, and perform the desired integration into elastic groups.

\subsection{Diffeological central extensions}
\label{sec:diff-centr-extens}

We recall the basic theory of central extensions of groups and Lie algebras. Fix a group $G$ and an abelian group $A$. We use concatenation to denote multiplication in $G$, and we use $+$ for the group addition in $A$. Implicitly we assume that $G$ acts trivially on $A$.

\begin{definition}
  A \define{central extension} of $G$ by $A$ consists of a group $\hat{G}$ and a short-exact sequence $A \hookrightarrow \hat{G} \twoheadrightarrow G$.  
\end{definition}

Central extensions of $G$ are classified by the second degree group cohomology. A 2-cochain on $G$ valued in $A$ is simply a map $G^2\to A$. The collection of 2-cochains forms a group with pointwise addition. The 2-cocycles and 2-coboundaries satisfy, respectively
\begin{align*}
  f(g_0g_1,g_2) + f(g_0,g_1) &=  f(g_0,g_1g_2) + f(g_1,g_2), \textrm{ and } \\
  f(g_0,g_1) &= -h(g_0g_1) + h(g_1) + h(g_0) \textrm{ for some } h\colon G \to A.
\end{align*}
Denote the resulting cohomology group by $H^2(G,A)$. We may insist that our cochains are normalized. Let
\begin{equation*}
  C^2(G,A) \coloneqq \{f\colon G^2 \to A \mid f(g_0,g_1) = 0 \textrm{ if } g_i = e \textrm{ for some } i\},
\end{equation*}
with subgroups of $n$-cocycles and $n$-coboundaries denoted by $Z^2(G,A)$ and $B^2(G,A)$, respectively. Then we still have $H^2(G,A) = Z^2(G,A)/B^2(G,A)$.

To each cocycle $f \in Z^2(G,A)$, we associate a central extension of $G$ by $A$, denoted $\hat{G}_f$ or $G \times_f A$. This has underlying set $G \times A$, and multiplication
\begin{equation*}
  (g,a)(h,b) \coloneqq (gh, a + b + f(g,h)).
\end{equation*}
The inclusion of $A$ into the second factor of $\hat{G}$, followed by the projection of the first factor to $G$, forms the required short-exact sequence. It is straightforward to check that $\hat{G}_f$ is a central extension. Moreover, $H^2(G,A)$ classifies central extensions of $G$ by $A$ up to equivalence.

The restriction to diffeological central extensions is fairly straightforward. Let $G$ be a diffeological group, and let $A$ be an abelian diffeological group.
\begin{definition}
  A \define{diffeological central extension} of $G$ by $A$ consists of a diffeological group $\hat{G}$ and a short-exact sequence $A \overset{\iota}{\hookrightarrow} \hat{G} \overset{\pi}{\twoheadrightarrow} G$, such that $\iota$ is an induction, and $\pi$ is a subduction.
\end{definition}
The conditions on $\iota$ and $\pi$ ensure that $A$ identifies with its image $\iota(A) \subseteq G$, viewed as a diffeological subgroup, and that $G$ identifies with $\hat{G}/\iota(A)$, where the latter carries the quotient diffeology. \emph{A fortiori} $\pi\colon \hat{G} \to G$ is a principal $A$-bundle. We caution that in $\cat{Dflg}$, principal bundles are not necessarily locally trivial, because $\pi$ does not necessarily admit a smooth local section.
\begin{example}
  \label{ex:7}
Fix some dense and countable subgroup $K$ of $\R$. Then $K \hookrightarrow \R \twoheadrightarrow \R/K$ expresses $\R$ as a central extension of $\R/K$ by $K$, where $K$ carries the subspace diffeology, and $\R/K$ carries the quotient diffeology. But the principal bundle $\R \to \R/K$ is not locally trivial.
\end{example}
This complicates the task of classifying diffeological group extensions by cohomology. On the other hand, we do not work in this generality.
\begin{definition}
  \label{def:central-extensions-2}
  Let $C_{s, \text{loc}}^2(G,A)$ denote those normalized group cochains which are smooth on a neighbourhood of $(e,\ldots, e) \in G^n$ of the form $U \times \cdots \times U$, for some $U \subseteq G$ open. This gives, suppressing the $(G,A)$ notation, cocycles $Z_{s, \text{loc}}^2$, coboundaries $B_{s, \text{loc}}^2$, and a cohomology $H_{s, \text{loc}}^2$.
\end{definition}
The ``s'' is for ``smooth,'' and ``loc'' is for ``local.''
\begin{remark}
The D-topology functor $D\colon \cat{Dflg} \to \cat{Top}$ does not generally preserve products; the natural map $D(X \times Y) \to D(X) \times D(Y)$ is only guaranteed to be continuous. This means that a neighbourhood of $(e, \ldots, e)$ in $D(G^n)$ need not contain any set of the form $U \times \cdots \times U$. Therefore the elements of $C_{s,\text{loc}}^2(G,A)$ are not merely maps $G^2 \to A$ that are smooth on some neighbourhood of $(e,e)$. Also observe that $D(G)$ need not be a topological group.
\end{remark}

For a connected diffeological group $G$, the cohomology $H^2_{s,\text{loc}}(G,A)$ classifies, up to equivalence, central extensions of $G$ by $A$ that are locally split, i.e.\ for which the projection $\pi$ admits a section $\sigma\colon G \to \hat{G}$ with $\sigma(e) = \hat{e}$ that is smooth in a neighbourhood of the identity. This is a consequence of two lemmas.

\begin{lemma}
\label{lem:28}
  Let $G$ be a group, and $U$ be a subset of $G$ equipped with a diffeology $\sr{D}_U$. Suppose we have $V \subseteq U$ a D-open subset containing the identity, such that
  \begin{enumerate}[label=(\roman*)]
  \item $V = V^{-1}$ and $V \cdot V \subseteq U$;
  \item the multiplication and inversion
    \begin{equation*}
      m\colon V \times V \to U, \quad i\colon V \to V 
    \end{equation*}
    are smooth (in the subspace diffeology of $V$, denoted $\sr{D}_V$);
  \item $V$ generates $G$ as a group.
  \end{enumerate}
  Then the finest group diffeology on $G$ containing $\sr{D}_V$ is
  \begin{equation*}
    \sr{D} \coloneqq \{p \mid p \textrm{ locally has the form } p = g \cdot q,\textrm{ for } \ g \in G, \ q \in \sr{D}_V\}.
  \end{equation*}
  Moreover
  \begin{enumeratea}
  \item The subset $V$ is a D-open subset of $(G, \sr{D})$.
    \item If $V' \subseteq U$ is another open neighbourhood of the identity satisfying (i) through (iii), then $\sr{D}$ is also the finest group diffeology containing $\sr{D}_{V'}$.
  \end{enumeratea}
\end{lemma}
We say that one diffeology is \define{finer} than another if it fewer plots. The proof of the lemma may seem long, but it employs standard techniques.
\begin{proof}
 The family $\{g \cdot q\}$ covers $G$. The plots of the diffeology generated by this family are those parametrizations $p$ which locally have the form $p = (g \cdot q) \circ F$, for some smooth map $F$ between open subsets of Cartesian spaces. But $q \circ F$ remains a plot of $\sr{D}_V$, so we see that the diffeology generated by $\{g \cdot q\}$ is precisely $\sr{D}$. This shows that $\sr{D}$ is a diffeology. Now we must show it is a group diffeology.
 
  We begin with the following observation. Take a plot $q$ of $V$, with $r_0 \in \dom q$. Then $q(r_0)^{-1} \cdot q$ is a plot of $U$ (because $m\colon V \times V \to U$ is smooth), sending $r_0$ to $e$. Since $V$ is a D-open neighbourhood of $e$, the restriction $q_1$ of $q(r_0)^{-1} \cdot q$ to some neighbourhood of $r_0$ is a plot of $V$.  Thus near $r_0$, we locally write $q = q(r_0) \cdot q_1$, for $q_1$ a plot of $V$ pointed at $e$.

  Now consider $g \in G$, and the map $q \cdot g$. Write $g = g_1\cdots g_l$, with $g_i \in V$, which is possible since $V$ generates $G$. Repeatedly applying the observation above, first to $q$ yielding $q_1$, then to $q_1 \cdot g_1$, yielding $q_2$, etc., we find that near $r_0 \in \dom q$,
  \begin{align*}
    q \cdot g &= q(r_0) \cdot (q_1 \cdot g_1) \cdots q_l \\
              &= q(r_0) \cdot g_1 \cdot (q_2 \cdot g_2) \cdots q_l \\
              &\cdots \\
              &= q(r_0) \cdot g \cdot q_{l+1}\\
              &\eqqcolon \ti{g} \cdot \ti{q}.
  \end{align*}
  The plot $\ti{q}$ is a plot of $V$ mapping $r_0$ to $e$.

We now check that the multiplication $m\colon G \times G \to G$ is smooth. It suffices to check that $g'q' \cdot gq$ is a plot of $\sr{D}$ for plots $q',q$ of $V$ on a shared domain, and $g',g \in G$. Fix $r_0 \in \dom q' = \dom q$. We may locally write $q' \cdot g$ as $\ti{g} \cdot \ti{q'}$, where $\ti{q'}$ is a plot of $V$ mapping $r_0$ to $e$. Then, locally
  \begin{equation*}
    g'q' \cdot gq = g' \ti{g} \cdot \ti{q'} \ti{q}, 
  \end{equation*}
  and $\ti{q'} q$ is a plot of $U$ since $m\colon V \times V \to U$ is smooth. Since $\ti{q'} q$ maps $r_0$ to $q(r_0) \in V$, this restricts to a plot of $V$, and we conclude that $g'q' \cdot gq$ is locally $g' \ti{g} \in G$ times a plot of $V$. Therefore the multiplication $m$ is smooth.

  We similarly show inversion is smooth. It suffices to show $(g \cdot q)^{-1}$ is a plot. Fix $r_0 \in \dom q$. We have $(g \cdot q)^{-1} = q^{-1} \cdot g^{-1}$, and $q^{-1}$ is a plot of $V$ because $i\colon V \to V$ is smooth. Near $r_0$, the product $q^{-1} \cdot g^{-1}$ therefore has the form $\ti{g} \cdot \ti{q^{-1}}$, for a plot $\ti{q^{-1}}$ of $V$. This shows $(g \cdot q)^{-1}$ is locally a plot of $\sr{D}$, hence is a plot.

  Thus we have verified that $\sr{D}$ is a group diffeology containing $\sr{D}_V$. Let $\sr{D}'$ denote the finest such diffeology. Then $\sr{D}' \subseteq \sr{D}$. Conversely, given $g \cdot q$ in the generating family for $\sr{D}$, we have $q \in \sr{D}_V \subseteq \sr{D}'$. Left multiplication $L_g\colon (G, \sr{D}') \to (G, \sr{D}')$ is smooth, so $g \cdot q = L_g \circ q$ is a plot of $\sr{D}'$. It follows that $\sr{D} \subseteq \sr{D}'$, as required.

  Let us now treat (a). To see $V$ is D-open, we may show that $(g \cdot q)^{-1}(V)$ is open for all $g \in G$ and $q \in \sr{D}_V$. If $g \in V$, denote by $L_g$ the map $(V, \sr{D}_V) \to (U, \sr{D}_U)$ given by multiplication by $g$. This is smooth. Then
  \begin{equation*}
    (g \cdot q)^{-1}(V) = \{r \mid q(r) \in (g^{-1} \cdot V) \cap V\} = \{r \mid q(r) \in L_g^{-1}(V)\},
  \end{equation*}
  and $L_g^{-1}(V)$ is D-open in $(V, \sr{D}_V)$ by smoothness of $L_g$. Therefore $(g \cdot q)^{-1}(V)$ is open. For general $g \in G$, at $r_0$ such that $g \cdot q(r_0) \in V$, write locally $q = q(r_0) \cdot q_1$, for a plot $q_1$ of $V$. Then $g \cdot q$ is, locally, $(g \cdot q(r_0)) \cdot q_1$, and we may apply the result of the previous case because $g \cdot q(r_0) \in V$. This provides an open neighbourhood of $r_0 \in \dom q$ which is sent by $g \cdot q$ into $V$, and thus proves that $(g \cdot q)^{-1}(V)$ is open in $\dom q$. Therefore $V$ is D-open.

  Finally, we treat (b). For $V' \subseteq U$ a D-open neighbourhood of the identity satisfying (i) through (iii), let $\sr{D}'$ denote the finest group diffeology containing $\sr{D}_{V'}$. Suppose first that $V' \subseteq V$. Then $\sr{D}' \subseteq \sr{D}$ is clear. Conversely, note that given $q \in \sr{D}_V$, and $r_0 \in \dom q$, we may locally write $q = q(r_0)^{-1} \cdot q_1$, where $q_1$ is a plot of $V'$. We may do so because $V'$ is D-open in $V$. This shows that $\sr{D} \subseteq \sr{D}'$. For general $V'$, consider $V \cap V'$. The finest diffeologies containing $\sr{D}_V$ and $\sr{D}_{V'}$, respectively, must coincide with the finest diffeology containing $V \cap V'$, hence with each other.
\end{proof}

\begin{lemma}
  \label{lem:central-extensions-3}
  Let $f \in Z_{s,\text{loc}}^2(G, A)$ be a locally smooth 2-cocycle. Then $\hat{G}_f = G \times_f A$ admits a group diffeology which makes it a locally split central extension of $G$ by $A$.
\end{lemma}

\begin{proof}
  Say $f\colon U \times U \to A$ is smooth, for an open neighbourhood $U$ of the identity. We may assume $U$ is connected. Choose an open neighbourhood $V$ of the identity in $U$ such that $V^{-1} = V$, and $V \cdot V \subseteq U$. Being an open neighbourhood of the identity in a topological group, $V$ generates $G$. Smoothness of $f$ on $V \times V$ implies that in $\hat{G}$,
  \begin{equation*}
    \textrm{multiplication}\colon  (V \times A) \times (V \times A) \to U \times A, \textrm{ and inversion}\colon (V \times A) \to (V \times A)
  \end{equation*}
  are smooth. Also, $V \times A$ generates $\hat{G}$. Then Lemma \ref{lem:28}, applied to the group $\hat{G}$, the diffeological space $U \times A$ (with the product diffeology), and the D-open subset $V \times A$, furnishes $\hat{G}$ with a group diffeology for which $V \times A$ is D-open.

  The inclusion $A \hookrightarrow \hat{G}$, given by $a \mapsto (e,a)$ is an induction into $V \times A$, hence into $\hat{G}$. Now we check that the projection $\pr\colon \hat{G} \to G$, given by $(g,a) \mapsto g$, is a subduction. Fix a plot $p$ of $G$, and take $r_0 \in \dom p$. Near $r_0$, we may write $p$ as $p = p(r_0)\cdot q$, for a plot $q$ of $V$. Then $(p(r_0),0)\cdot (q,0)$ is the required plot of $\hat{G}$ which projects to $p(r_0) q = p$.

   Finally, the map $V \to V \times A$ given by $g \mapsto (g,0)$ is smooth, thus it is smooth as a map $V \to \hat{G}$. This is a section of the projection, and its existence corresponds to a local trivialization of $\hat{G} \to G$. Note that this section is not a group homomorphism, so it may not extend to a smooth map beyond $V$.
 \end{proof}

 \begin{proposition}
   \label{prop:central-extensions-3}
   Let $G$ be a connected diffeological group, and $A$ an abelian diffeological group. The cohomology $H^2_{s,\text{loc}}(G,A)$ classifies, up to equivalence, the locally split central extensions of $G$ by $\hat{G}$.
 \end{proposition}
 \begin{proof}
   Given $f \in Z_{s,\text{loc}}^2(G,A)$, Lemma \ref{lem:central-extensions-3} gives the corresponding locally split central extension of $G$ by $A$. Conversely, given such an extension, let $\sigma\colon G \to A$ be a section of the projection that is smooth near $e$. Then
   \begin{equation*}
     f(g,h) \coloneqq \sigma(g)\sigma(h)\sigma(gh)^{-1}
   \end{equation*}
   is the element of $Z_{s,\text{loc}}^2(G,A)$ that corresponds to the original $\hat{G}$.
 \end{proof}

Now we turn to Lie algebras. Let $\fk{g}$ be a Lie algebra, and $\fk{a}$ be a trivial Lie algebra, thus merely a vector space. We implicitly assume that $\fk{g}$ acts trivially on $\fk{a}$.
\begin{definition}
A \define{central extension} of $\fk{g}$ by $\fk{a}$ consists of a Lie algebra $\hat{\fk{g}}$ and a short-exact sequence of Lie algebras $\fk{a} \hookrightarrow \hat{\fk{g}} \twoheadrightarrow \fk{g}$.  
\end{definition}
Like with groups, central extensions of Lie algebras are classified by second-degree cohomology. A 2-cocycle for $\fk{g}$ valued in $\fk{a}$ is an antisymmetric bilinear map $\omega\colon \fk{g} \times \fk{g} \to \fk{a}$ satisfying $\sum_{\circlearrowright} \omega([x,y],z) = 0$.\footnote{We use $\sum_\circlearrowright$ to denote a sum over cyclic permutations of the arguments.} A 2-coboundary is a 2-cocycle which, additionally, is of the form $\omega(x,y) = b([x,y])$, for some linear map $b\colon \fk{g} \to \fk{a}$. We use $Z^2(\fk{g}, \fk{a})$, $B^2(\fk{g}, \fk{a})$, and $H^2(\fk{g}, \fk{a})$, to denote the 2-cocycles, 2-coboundaries, and the resulting cohomology group, respectively.

  To each $\omega \in Z^2(\fk{g}, \fk{a})$, we associate a central extension of $\fk{g}$ by $\fk{a}$, denoted $\hat{\fk{g}}_\omega$ or $\fk{g} \oplus_\omega \fk{a}$. This has underlying vector space $\fk{g} \oplus \fk{a}$, and Lie bracket
  \begin{equation*}
  {[(x,a),(y,b)]} \coloneqq ({[x,y]}, \omega(a,b)).
\end{equation*}
The inclusion of $\fk{a}$ into the second factor of $\hat{\fk{g}}_\omega$, followed by the projection of the first factor to $\fk{g}$, forms the required short exact sequence. Like with groups, the second degree Lie algebra cohomology $H^2(\fk{g}, \fk{a})$ classifies central extensions of $\fk{g}$ by $\fk{a}$ up to equivalence.

We now let $\fk{g}$ be a diffeological Lie algebra, and $\fk{a}$ be a trivial diffeological Lie algebra.

\begin{definition}
  A \define{diffeological central extension} of $\fk{g}$ by $\fk{a}$ is a Lie algebra central extension of $\fk{g}$ by $\fk{a}$, which is simultaneously a diffeological central extension of the underlying abelian diffeological groups.
\end{definition}

Since central extensions of diffeological Lie algebras are also central extensions of diffeological groups, the same considerations apply. However, we will be content to work with globally split central extensions, and thus do not need versions of Lemmas \ref{lem:28} or \ref{lem:central-extensions-3}. We have, suppressing $(\fk{g}, \fk{a})$ from the notation, the groups $Z_s^2$ and $B_s^2$ of smooth 2-cocycles and 2-coboundaries, and this gives the smooth Lie algebra cohomology $H_s^2$. Given $\omega \in Z_s^2$, the associated central extension $\hat{\fk{g}}_\omega$, equipped with the product diffeology, is a split diffeological central extension of $\fk{g}$ with values in $\fk{a}$. The cohomology $H_s^2$ classifies such central extensions up to equivalence.

\begin{proposition}
  \label{prop:23}
Let $\fk{g}$ be a diffeological Lie algebra, and $\fk{a}$ be a trivial diffeological Lie algebra.  The second degree smooth Lie algebra cohomology $H_s^2(\fk{g}, \fk{a})$ classifies split diffeological central extensions of $\fk{g}$ by $\fk{a}$, up to equivalence.
\end{proposition}

The following table summarizes our discussion.
\begin{table}[h]
  \centering
  \begin{tabular}[h]{l | l | l}
    Extension & Cohomology & Correspondence \\ \hline
    Central extensions of groups & $H^2(G,A)$ & $f \mapsto \hat{G}_f \coloneqq G \fiber{}{f} A$ \\
    \makecell[l]{Locally split central extensions \\ \quad of diffeological groups} & $H^2_{s, \text{loc}}(G,A)$ & $f \mapsto \hat{G}_f$, diffeology in Lemma \ref{lem:central-extensions-3} \\
    Central extensions of Lie algebras & $H^2(\fk{g},\fk{a})$ & $\omega \mapsto \hat{\fk{g}}_\omega \coloneqq \fk{g} \oplus_\omega \fk{a}$ \\
    \makecell[l]{Split central extensions \\ \quad of diffeological Lie algebras} & $H^2_s(\fk{g},\fk{a})$ & $\omega \mapsto \hat{\fk{g}}_\omega$, product diffeology.
  \end{tabular}
\end{table}

To end this section, we differentiate locally smooth cocycles. Let $G$ be an elastic group, with Lie algebra $\fk{g}$. Our differentiation of a 2-cocycle $f \in Z_{s, \text{loc}}^2(G,A)$ to a 2-cocycle $L(f) \in Z_s^2(\fk{g}, \fk{a})$ amounts to taking the Hessian of $f$ at $(e,e)$. Suppose that $M, N$, and $P$ are elastic spaces, and $f\colon M \times N \to P$ is smooth. Then we have the map
\begin{equation*}
  T_{(1)}T_{(2)}f\colon TM \times TN \to T^2P.
\end{equation*}
If, moreover, we have $m_0 \in M$, and $n_0 \in N$, and $p_0 \in P$ such that $f(m_0,\cdot) = f(\cdot, n_0) = p_0$, then Lemma \ref{lem:appendix-3} implies that, whenever $v \in T_{m_0}M$ and $w \in T_{n_0}N$,
\begin{equation*}
  T_{(1)}T_{(2)}f(v,w) = T_{(2)}T_{(1)}f(v,w) \in \lambda(T_{p_0}P) \cong T_{p_0}P.
\end{equation*}
The map $T_{(1)}T_{(2)}f$ is also bilinear as a map $T_{m_0}M \times T_{n_0}N \to T_{p_0}P$.

\begin{definition}
  Let $f \in Z_{s, \text{loc}}^2(G,A)$ be a 2-cocycle. We denote
  \begin{equation*}
    d^2f \coloneqq T_{(1)}T_{(2)}f \text{ restricted to a map } \fk{g} \times \fk{g} \to \fk{a},
  \end{equation*}
  and we define $L(f)\colon \fk{g} \times \fk{g} \to \fk{a}$ by
  \begin{equation*}
    L(f)(\xi, \eta) \coloneqq T_{(1)}T_{(2)}f(\xi, \eta) - T_{(2)}T_{(1)}f(\eta, \xi).
  \end{equation*}
\end{definition}

Now we can differentiate cocycles.
\begin{proposition}
  \label{prop:19}
Let $G$ be an elastic diffeological group, and $A$ be an abelian elastic diffeological group. Let $f \in Z^2_{s,\text{loc}}(G,A)$. Then $L(f)$ is a smooth Lie algebra 2-cocycle for $\fk{g}$ with values in $\fk{a}$, and the Lie algebra of $\hat{G}_f$ is $\hat{\fk{g}}_{L(f)}$.
\end{proposition}

\begin{proof}
  Suppose that $f$ is smooth on $U \times U$, for $U$ an open neighbourhood of the identity in $G$. Since $d^2f$ is a smooth bilinear map, so is $L(f)$. Anti-symmetry is clear. Now we evaluate the Lie bracket on $\hat{G}_f$. Recall that the conjugation in $\hat{G}_f$ is explicitly given by
  \begin{equation*}
    (g,z) (h, w) (g, z)^{-1} = (ghg^{-1}, w + f(g,h) - f(ghg^{-1},g)).
  \end{equation*}
  Applying $T_{(1)}T_{(2)}$ and evaluating at $((\xi, V), (\eta, W))$ yields, on the left side, $[(\xi, V),(\eta, W)]$ by Lemma \ref{lem:appendix-2}. On the right side, we have $[\xi, \eta]$ in the first slot, and the $f(g,h)$ term gives $T_{(1)}T_{(2)}f(\xi, \eta)$. For the $-f(ghg^{-1},g)$ term, we consider
  \begin{equation*}
    f \circ (c, \pr_1)\colon G \times G \to G.
  \end{equation*}
  Here $c$ denotes conjugation. Applying $T_{(1)}T_{(2)}$ to this gives
  \begin{align*}
    T_{(1)}T_{(2)}(f \circ (c, \pr_1)) &= T_{(1)}(Tf \circ (Tc,\pr_1) \circ (0 \times 1)) \\
                                       &= T_{(1)}(Tf \circ (T_{(2)}c, 0 \circ \pr_1)) \\
                                       &= T_{(1)}(T_{(1)}f \circ (T_{(2)}c, \pr_1)) \\
                                       &= T(T_{(1)}f \circ (T_{(2)}c,\pr_1)) \circ (1 \times 0) \\
    &= TT_{(1)}f \circ (T_{(1)}T_{(2)}c, \pr_1).
  \end{align*}
  Apply Lemma \ref{lem:appendix-1} to break this into a sum:
  \begin{align*}
    TT_{(1)}f \circ (T_{(1)}T_{(2)}c, \pr_1) &= T_{(1)}T_{(1)}f \circ (T_{(1)}T_{(2)}c, \pi \circ \pr_1) \\
    &\qquad + T_{(2)}T_{(1)}f \circ (\pi \circ  T_{(1)}T_{(2)}c, \pr_1).
  \end{align*}
  The first summand will vanish, because $f(e,\cdot) = f(\cdot, e) = e$, and the second is $T_{(2)}T_{(1)}f \circ (\pr_2, \pr_1)$. Thus upon differentiating, the $-f(ghg^{-1},g)$ term gives $T_{(2)}T_{(1)}f(\eta, \xi) = T_{(1)}T_{(2)}f(\eta, \xi)$. Therefore
  \begin{align*}
    [(\xi, V), (\eta, W)] &= ([\xi, \eta], T_{(1)}T_{(2)}f(\xi, \eta) - T_{(1)}T_{(2)}f(\eta,\xi)) \\
    &= ([\xi,\eta], L(f)(\xi, \eta)).
  \end{align*}
  The fact that the Lie bracket can be expressed in this way implies that $L(f)$ is a Lie algebra 2-cocycle in $Z^2_s(\fk{g}, \fk{a})$, and that $\Lie(\hat{G}_f) = \hat{\fk{g}}_{L(f)}$.
\end{proof}

\subsection{Integrating convenient Lie algebras}
\label{sec:integr-conv-lie}

Having established how to differentiate locally smooth 2-cocycles, we turn to the problem of integration. We draw from results in the literature \cite{Neeb02,Woc11}. These authors work in the category of infinite-dimensional manifolds modelled on lctvs, and Michel-Bastiani smooth maps between them (cf.\ Remark \ref{rem:convenient-2}), but we will apply their results to convenient manifolds.

For our preliminary discussion, fix a connected and simply connected convenient Lie group $G$, with Lie algebra $\fk{g}$, and a convenient abelian Lie algebra $\fk{a}$. Also take a 2-cocycle $\omega \in Z^2_s(\fk{g}, \fk{a})$. Let $\omega^l \in \Omega^2(G, \fk{a})$ be the $\fk{a}$-valued, $G$-invariant differential 2-form with $\omega^l(e) = \omega$.\footnote{For a diffeological vector space $V$, and a diffeological space $X$, we have not defined the meaning of a $V$-valued differential form on $X$. Instead, to understand $\omega^l$ we simply remain in the setting of convenient manifolds.}

\begin{definition}
The \define{period homomorphism} is the map
\begin{equation*}
  \operatorname{per}_\omega\colon H_2(G) \to \fk{a}, \quad [\sigma] \mapsto \int_\sigma \omega^l.
\end{equation*}
We denote its image by $\Pi_\omega \coloneqq \operatorname{per}_\omega(H_2(G))$, and we have the quotient $\pi\colon \fk{a} \to \fk{a}/\Pi_\omega$. We set $A \coloneqq \fk{a}/\Pi_\omega$.
\end{definition}

Neeb, inspired by van Est \cite{Est58}, gave the following construction of a locally smooth 2-cocycle integrating $\omega$. While Neeb assumes that $\Pi_\omega$ is a topologically discrete subgroup of $\fk{a}$, we make the less-restrictive assumption that $\Pi_\omega$ is a diffeologically discrete subgroup of $\fk{a}$.

\begin{lemma}[{\cite[Lemmas 6.2,6.3]{Neeb02}}]
  \label{lem:central-extensions-2}
  Let $u \subseteq \fk{g}$ be an open convex neighbourhood of the origin, with a chart $\varphi\colon u \to G$ sending $0$ to $e$ and such that $T_0\varphi = \id_{\fk{g}}$. Define, for $x \in u$,
  \begin{equation*}
    \alpha_{\varphi(x)}\colon [0,1] \to G, \quad \alpha_{\varphi(x)}(t) \coloneqq \varphi(tx).
  \end{equation*}
  Let $v \subseteq u$ be an open convex neighbourhood of the origin, such that $\varphi(v)\varphi(v) \subseteq \varphi(u)$, and define, for $x,y \in V$,
  \begin{equation*}
  x * y \coloneqq \varphi^{-1}(\varphi(x)\varphi(y)).  
  \end{equation*}
  We use this multiplication to define the map $\gamma_{x,y}$ from the standard 2-simplex $\Delta^2$
  \begin{equation*}
    \gamma_{x,y}\colon \Delta^2 \to u, \quad \gamma_{x,y}(t,s) \mapsto t(x*sy) + s(x*(1-t)y),
  \end{equation*}
  and we define $\sigma_{x,y} \coloneqq \varphi \circ \gamma_{x,y} \colon \Delta^2 \to G$. 
  \begin{enumeratea}
  \item The function
    \begin{equation*}
      f_0\colon v \times v \to \fk{a}, \quad f_0(x,y) \coloneqq \int_{\sigma_{x,y}}\omega^l
    \end{equation*}
    is smooth, and $d^2f_0 = \frac{1}{2}\omega$.
  \item The function
    \begin{equation*}
      f\colon \varphi(v) \times \varphi(v) \to A, \quad f(\varphi(x),\varphi(y)) \coloneqq  \pi(f_0(x,y))
    \end{equation*}
    is a locally smooth 2-cocycle in $Z^2_{s, \text{loc}}(G,A)$, and $L(f) = \omega$.
  \end{enumeratea}
\end{lemma}

\begin{proof}
  Part (a) is simply a restatement of \cite[Lemma 6.2]{Neeb02}, under identifying $T_0v \cong \fk{g}$ and $T_0\fk{a} \cong \fk{a}$. For (b), we begin with the fact that, since $T_0\varphi = \id_{\fk{g}}$,
  \begin{equation*}
   L(f_0 \circ (\varphi^{-1}, \varphi^{-1})) = \omega.
 \end{equation*}
Now, by Corollary \ref{cor:groupoids-1} the map $T_0\pi\colon \Lie(\fk{a}) \cong \fk{a} \to \Lie(A)$ is an isomorphism. Therefore $L(f) = \omega$ also holds.
\end{proof}

We now arrive at our main result for this section.

\begin{maintheorem}{III}
  \label{thm:central-extensions-1}
  Suppose that $\fk{a} \hookrightarrow \hat{\fk{g}} \twoheadrightarrow \fk{g}$ is a topologically split central extension of convenient Lie algebras. Fix $\omega \in Z_s^2(\fk{g}, \fk{a})$ corresponding to this extension. Assume that there is some connected and simply-connected convenient Lie group $G$ with $\Lie(G) = \fk{g}$. If the image $\Pi_\omega$ of the associated period homomorphism $\per_\omega\colon H_2(G) \to \fk{a}$ is diffeologically discrete, then there is some elastic diffeological group $\hat{G}$ centrally extending $G$ by $\fk{a}/\Pi_\omega$, for which
  \begin{equation*}
    \begin{tikzcd}
      \fk{a}/\Pi_\omega \ar[d, dashed, "\Lie"] \ar[r, hook] & \hat{G} \ar[r, two heads] \ar[d, dashed, "\Lie"] & G \ar[d, dashed, "\Lie"]  \\
      \fk{a} \ar[r, hook] & \hat{\fk{g}} \ar[r, two heads] & \fk{g}.
    \end{tikzcd}
  \end{equation*}
\end{maintheorem}
\begin{proof}
Let $A \coloneqq \fk{a}/\Pi_\omega$.  Appealing to Lemma \ref{lem:central-extensions-2}, we get $f \in Z_{s,\text{loc}}^2(G,A)$ with $L(f) = \omega$. Then, by Proposition \ref{prop:19}, the central extension of $G$ by $A$ corresponding to $f$ differentiates to the central extension of $\fk{g}$ by $\fk{a}$ corresponding to $\omega$. The former gives $\hat{G} \coloneqq \hat{G}_f$; the latter is $\hat{g} \cong \hat{g}_\omega$.
\end{proof}

With this machinery, we integrate a class of convenient Lie algebras, as promised. Let $\fk{z}$ denote the center of $\fk{g}$, and set $\fk{g}_{\text{ad}} \coloneqq \fk{g}/\fk{z}$. If $\fk{z} \hookrightarrow \fk{g} \twoheadrightarrow \fk{g}_{\text{ad}}$ is a split central extension of Lie algebras, it corresponds to the cocycle
\begin{equation*}
  \omega\colon \fk{g}_{\text{ad}} \times \fk{g}_{\text{ad}} \to \fk{z}, \quad (\xi, \eta) \mapsto [\xi, \eta]_{\fk{z}},
\end{equation*}
where $[\cdot,\cdot]_{\fk{z}}$ denotes the $\fk{z}$-component of the Lie bracket. If there is, furthermore, a connected and simply-connected convenient Lie group $G_{\text{ad}}$ with $\Lie(G_{\text{ad}}) = \fk{g}_{\text{ad}}$, then we have the associated period homomorphism $\per_\omega\colon H_2(G_{\text{ad}}) \to \fk{z}$. If this has diffeologically discrete image, then Theorem \ref{thm:central-extensions-1} gives $G$ with $\Lie(G) = \fk{g}$.

\begin{remark}
  We require the assumption that the central extension $\fk{z} \hookrightarrow \fk{g} \twoheadrightarrow \fk{g}_{\text{ad}}$ is split in order to identify the extension with a 2-cocycle in $Z_s^2(\fk{g}_{\text{ad}}, \fk{z})$. If the extension is not split, Neeb \cite{Neeb06} has given a two-step process that salvages the strategy. One first integrates the central term in a split-exact sequence
  \begin{equation*}
    \fk{g} \hookrightarrow \hat{\fk{g}} \coloneqq \fk{g} \oplus_b \fk{g}_{\text{ad}} \twoheadrightarrow |\fk{g}_{\text{ad}}|,
  \end{equation*}
  where $|\cdot|$ indicates the underlying vector space and $b \in Z_s^2(\fk{g}, \fk{g}_{\text{ad}})$ is the 2-cocycle induced by the bracket. Say this integration is $\hat{G}$. Then, using the fact that $\fk{g}$ is a closed ideal of $\hat{\fk{g}}$, argue that the inclusion $\fk{g} \hookrightarrow \hat{\fk{g}}$ lifts to an inclusion $G \hookrightarrow \hat{G}$ for some $G$.

  Thus, removing our assumption that the central extension is split amounts to understanding when, given an elastic group $G$ and a closed ideal $\fk{h}$ of $\Lie(G)$, there is some subgroup $H$ of $G$ with $\Lie(H) = \fk{h}$. 
\end{remark}

The various pre-conditions of Theorem \ref{thm:central-extensions-1}, applied to $\fk{z} \hookrightarrow \fk{g} \twoheadrightarrow \fk{g}_{\text{ad}}$, hold in the following contexts:
\begin{itemize}
\item If $\fk{g}$ is a Banach Lie algebra, and $\fk{z}$ is finite-dimensional, then $\fk{z} \hookrightarrow \fk{g} \twoheadrightarrow \fk{g}_{\text{ad}}$ is split. This follows from the facts that, by the Hahn-Banach theorem, $\fk{z}$ has a closed complement, and by the open mapping theorem, the existence of this closed complement is equivalent to the existence of a continuous (hence smooth) section $\fk{g}_{\text{ad}} \to \fk{g}$.
\item If $\Pi_\omega$ is countable, then it is diffeologically discrete.
\item If $\fk{g}_{\text{ad}}$ is a ``locally exponential'' convenient Lie algebra, then Neeb \cite[Theorem IV.3.8]{Neeb06}\footnote{note that this citation contains no proof} has constructed a convenient Lie group integrating it, and we may then pass to the simply-connected cover. In particular, since Banach Lie algebras are locally exponential, when $\fk{g}$ is a Banach Lie algebra we may always find a connected and simply-connected Banach Lie group $G_{\text{ad}}$ with $\Lie(G_{\text{ad}}) = \fk{g}_{\text{ad}}$.  
\end{itemize}

We finish with the classic examples of non-enlargeable Banach Lie algebras, given by van Est and Korthagen \cite{EstKor64} and Douady and Lazard \cite{DouadLaz66}, as retold by Neeb \cite{Neeb06}.
\begin{example}[{\cite[Section 6]{EstKor64}}]
  Consider the Banach-Lie algebra $C^1(S^1, \fk{su}_2(\C))$, and the $\R$-valued cocycle
  \begin{equation*}
    \omega(f,g) \coloneqq \int_0^1 \operatorname{trace}(f(t)g'(t))dt.
  \end{equation*}
  For the integral, we identify the continuous maps from $S^1$ with maps from $\R$ with period $1$. Let $\fk{g}_{\text{EK}}$ denote the central extension of $C^1(S^1,\fk{su}_2(\C))$ by $\omega$. One verifies that
  \begin{equation*}
    (\fk{g}_{\text{EK}})_{\text{ad}} \cong C^1(S^1, \fk{su}_2(\C)) \text{ and } G_{\text{ad}} \cong C^1(S^1, SU_2(\C)),
  \end{equation*}
  so then $\pi_2(G_{\text{ad}}) = \pi_3(SU_2(\C)) = \pi_3(S^3) = \Z$. This means that $\Pi_\omega$ is at most countable. Therefore Theorem \ref{thm:central-extensions-1} applies, integrating $\fk{g}_{\text{EK}}$ to an elastic group $G_{\text{EK}}$.
\end{example}

\begin{example}[{\cite[Section IV]{DouadLaz66}}]
  Let $H$ be a separable infinite-dimensional complex Hilbert space. Consider the Banach-Lie algebra $\fk{u}(H)$ of skew-adjoint bounded linear operators on $H$. This has center $i\R$. For an irrational number $\theta$, define
  \begin{equation*}
    \fk{n}_\theta \coloneqq (1_H, \theta \cdot 1_H) \cdot i\R.
  \end{equation*}
  This is a closed ideal of the center $i\R \times i\R$ of $\fk{u}(H) \times \fk{u}(H)$, thus $\fk{g}_{\text{DL}} \coloneqq (\fk{u}(H) \times \fk{u}(H))/\fk{n}_\theta$ is a Banach-Lie algebra. The central extension
  \begin{equation*}
    i\R \cong \fk{z}_{\text{DL}} \hookrightarrow \fk{g}_{\text{DL}} \twoheadrightarrow (\fk{g}_{\text{DL}})_{\text{ad}} \cong (\fk{u}(H) \times \fk{u}(H))_{\text{ad}}
  \end{equation*}
  is split because $\fk{z}$ is finite-dimensional. Furthermore, by the discussion in \cite[Pages 415--416]{Neeb06}, the corresponding period group $\Pi_\omega$ is isomorphic to $\Z^2$, hence is countable. Therefore, again, we may use Theorem \ref{thm:central-extensions-1} to integrate $\fk{g}_{\text{DL}}$ to an elastic group $G_{\text{DL}}$.
\end{example}
\begin{remark}
  In \cite{BelPel25}, Belti\c{t}\u{a} and Pelletier generalize Banach Q-manifolds to H-manifolds. Both are sets equipped with atlases of some kind. Belti\c{t}\u{a} and Pelletier integrate every Banach-Lie algebra to an H-group, and remark \cite[Example 4.9]{BelPel25} that it is unclear whether Douady and Lazard's example $\fk{g}_{\text{DL}}$ integrates to a Q-group. However, careful examination of our integration, specifically Lemma \ref{lem:central-extensions-3}, reveals that $G_{\text{DL}}$ is locally diffeomorphic to $(\fk{z}_{\text{DL}}/\Pi_\omega) \times V$, where $V$ is an open subset of a Banach-Lie group. Since $\Pi_\omega$ is a countable subgroup of $\fk{z}_{\text{DL}}$, the set $(\fk{z}_{\text{DL}}/\Pi_\omega) \times V$, hence the set $G_{\text{DL}}$, also inherit the structure of a Q-group. We anticipate that this Q-group differentiates to $\fk{g}_{\text{DL}}$ under Belti\c{t}\u{a} and Pelletier's Lie functor.

  There are natural faithful functors from the categories of Q- and H-manifolds into $\cat{Dflg}$. Knowing whether these functors are full, and having a characterization of their essential images, would be helpful in comparing Belti\c{t}\u{a} and Pelletier's approach to the diffeological one.
\end{remark}

\appendix

\section{Tangent categories}
\label{sec:tangent-categories}

In this Appendix, we collect some facts about tangent categories. Specifically, Lemmas \ref{lem:appendix-3} and \ref{lem:appendix-2} hold in an arbitrary tangent category, so we prove them in that generality. We include the definitions and axioms of a tangent category for completeness rather than pedagogy. We direct the reader to \cite{AinBloh25} for a gentler introduction to tangent categories. The articles \cite{CocCrut14} and \cite{CocCrut15} are also good sources.

\subsection{First definitions}
\label{sec:first-definitions}

Let $\cat{C}$ be a category, with a terminal object $*$. We assume that the overcategory $\cat{C}\downarrow X$ (hence also $\cat{C} \cong \cat{C} \downarrow *$) has all finite products. We need three notions to define a tangent category.
\begin{definition}
  A \define{bundle of abelian groups} $\pi\colon A \to X$ is an abelian group object in $\cat{C}\downarrow X$. 
\end{definition}
The data of such a bundle consists of the bundle projection $\pi\colon A \to X$, the addition $+\colon A \fiber{}{X} A \to A$, zero-section $0\colon A \to A$, and inversion $i\colon A \to A$.

\begin{definition}
  Let $T\colon \cat{C} \to \cat{C}$ be an endofunctor, and $\tau\colon T^2 \to T^2$ a natural transformation. Let $\tau_{12} \coloneqq \tau T$ and $\tau_{23} \coloneqq T \tau$ be the trivial extensions of $\tau$ to natural transformations $T^3 \to T^3$. We call $\tau$ a \define{symmetric structure} if $\tau \circ \tau = 1$, and $\tau_{12} \circ \tau_{23} \circ \tau_{12} = \tau_{23} \circ \tau_{12} \circ \tau_{23}$.
\end{definition}

\begin{definition}
  We say an endofunctor $T\colon \cat{C} \to \cat{C}$ \define{preserves fiber products} of a bundle of abelian groups $p\colon A \to X$ if for all $k \geq 1$, the natural morphisms
  \begin{equation*}
    \nu_{k,X} \colon T(A \fiber{}{X} \cdots \fiber{}{X} A) \to TA \fiber{}{TX} \cdots \fiber{}{TX} TA
  \end{equation*}
  are all isomorphisms.
\end{definition}

Here is the definition of a tangent structure.
\begin{definition}
  The data of a \define{tangent structure} (with negatives) on $\cat{C}$ consists of an endofunctor $T\colon \cat{C} \to \cat{C}$, and natural transformations
  \makeatletter
  \@fleqntrue
  \begin{alignat*}{2}
    & \textrm{the \define{footprint}}\quad &  \pi &\colon T \to 1 \\
    &\textrm{the \define{zero-section}}\quad & 0 &\colon 1 \to T \\
    & \textrm{the \define{addition}}\quad & +&\colon T_2 \to T \\
    &\textrm{the \define{vertical lift}}\quad & \lambda&\colon T \to T^2 \\
    &\textrm{the \define{symmetric structure}}\quad & \tau &\colon T^2 \to T^2.  
  \end{alignat*}
  \makeatother
  The axioms are:
  \begin{itemize}
  \item The footprint $\pi\colon T \to 1$ is a bundle of abelian groups over $1$, with zero-section $0$ and addition $+$.
  \item All pullbacks $T_k \coloneqq T \fiber{}{1} \cdots \fiber{}{1} T$ exists, and are pointwise preserved by $T$.
  \item The transformation $\tau$ is both a symmetric structure and a morphism of bundles of abelian groups.
  \item The diagrams
    \begin{equation*}
      \begin{tikzcd}
        T  & T^2 \\
        1  & T
        \ar["\lambda", from=1-1, to=1-2]
        \ar["\pi"', from=1-1, to=2-1]
        \ar["\pi T", from=1-2, to=2-2]
        \ar["0", from=2-1, to=2-2]
      \end{tikzcd}
      \quad
      \begin{tikzcd}
        T & T^2   \\
        T^2  & T^3
        \ar["\lambda", from=1-1, to=1-2]
        \ar["\lambda"', from=1-1, to=2-1]
        \ar["\lambda T", from=1-2, to=2-2]
        \ar["T\lambda", from=2-1, to=2-2]
      \end{tikzcd}
    \end{equation*}
    commute, and the first is a morphism of abelian groups.
  \item The diagrams
    \begin{equation*}
      \begin{tikzcd}
        & T & \\
        T^2  & & T^2
        \ar["\lambda"', from=1-2, to=2-1] \\
        \ar["\lambda", from=1-2, to=2-3]
        \ar["\tau", from=2-1, to=2-3]
      \end{tikzcd}
      \quad
      \begin{tikzcd}
        T^2  & T^3  & T^3 \\
        T^2  & & T^3
        \ar["T\lambda", from=1-1, to=1-2]
        \ar["\tau T", from=1-2, to=1-3]
        \ar["\tau"', from=1-1, to=2-1]
        \ar["\lambda T", from=2-1, to=2-3]
        \ar["T\tau", from=1-3, to=2-3]
      \end{tikzcd}
    \end{equation*}
    commute.
  \item The diagram
    \begin{equation*}
      \begin{tikzcd}
        T_2  & T_2T & T^2 & T^2 \\
        1  & & & T
        \ar["T0 \fiber{}{0} \lambda", from=1-1, to=1-2]
        \ar["+_T", from=1-2, to=1-3]
        \ar["\tau", from=1-3, to=1-4]
        \ar["\pi \circ \pr_1"', from=1-1, to=2-1]
        \ar["T\pi", from=1-4, to=2-4]
        \ar["0", from=2-1, to=2-4]
      \end{tikzcd}
    \end{equation*}
    is a pointwise pullback. We denote the top transformation by $\lambda_2$.
  \end{itemize}    
\end{definition}
\begin{remark}
  The last axiom, in the presence of all others, is by \cite[Lemma 3.10]{CocCrut14}, equivalent to the easier-to-write condition that    
  \begin{equation*}
    \begin{tikzcd}
      T & T^2  \\
      1  & T_2
      \ar["\lambda", from=1-1, to=1-2]
      \ar["\pi"', from=1-1, to=2-1]
      \ar["{(\pi T, T\pi)}", from=1-2, to=2-2]
      \ar["{(0,0)}", from=2-1, to=2-2]
    \end{tikzcd}
  \end{equation*}
  is a pointwise pullback. However, it is $\lambda_2$ that will appear in the definition of the Lie bracket.
\end{remark}

The tangent structure lets us define a Lie bracket. Fix $v$ and $w$ vector fields on $X$, i.e.\ sections of $\pi_X\colon TX \to X$. The footprint in $TX$ of $Tw \circ v$ is $w$:
\begin{equation*}
  \pi_{TX} \circ Tw \circ v = w \circ \pi_X \circ v = w \circ 1_X = w.
\end{equation*}
Using the fact $\pi_{TX} \circ \tau_X = T\pi_X$, we see the footprint of $\tau_X \circ Tv \circ w$ is $w$:
\begin{equation*}
  \pi_{TX} \circ \tau_X \circ Tv \circ w = T\pi_X \circ Tv \circ w = T1_X \circ w = 1_{TX} \circ w = w.
\end{equation*}
Thus we may subtract them in $T^2X$. In other words, we have a map
\begin{equation*}
  \delta(v,w)\colon X \to T^2X, \quad \delta(v,w) = Tw \circ v - \tau_X \circ Tv \circ w.
\end{equation*}
On the other hand, again using the chain rule, we can check that
\begin{equation*}
  T\pi_X \circ \delta(v,w) = v-v = 0_X.
\end{equation*}
Thus $\delta(v,w)$ takes its values in the kernel of $T\pi_X$, which, by the axioms of a tangent category, coincides with the image of $\lambda_{2,X}$.
\begin{definition}
  \label{def:tangent-structure-2}
The Lie bracket of vector fields $v$ and $w$ in $\fk{X}(X)$ is defined as the unique vector field $[v,w]$ such that
  \begin{equation*}
    \delta(v,w) = \lambda_{2,X} \circ (w, [v,w]).
  \end{equation*}
\end{definition}

Rosick\'{y} \cite{Ros84} announced the result that the bracket defined above satisfies the Jacobi identity. Cockett and Crutwell \cite{CocCrut15}, working from Rosick\'{y}'s notes, supplied a complete proof.

\begin{theorem}[Rosick\'{y}; Cockett-Crutwell]
  \label{thm:appendix-1}
For an object $X$ of a tangent category $\cl{C}$, the Lie bracket on $\fk{X}(X)$ satisfies the Jacobi identity.  
\end{theorem}

We now add more structure to our tangent category. Recall that $\cat{C}$ has finite products.
\begin{definition}
  The tangent category $(\cat{C}, T)$ is \define{Cartesian} if $T* \cong *$ and the natural morphisms $T(X\times Y) \to TX \times TY$ are isomorphisms for all objects $X$ and $Y$.
\end{definition}

Fix a ring object $R$ in a Cartesian tangent category $\cat{C}$. The natural transformation $1 \times R \to 1$ inherits from $R$ the structure of a ring object in the category $\operatorname{End}(\cat{C}) \downarrow 1$.

\begin{definition}
  An $R$-\define{module structure on }$T$ is a $(1 \times R \to 1)$-module structure on $\pi$.
\end{definition}
Such a structure is given by a natural morphism $\kappa\colon R \times T \to T$, and supplies the $R$-multiplication on the fiber. We have a morphism $TR \to R \times R$. If this is an isomorphism of $R$-modules, we call $R$ \define{tangent-stable}.
\begin{definition}[{\cite[Definition 2.24]{AinBloh25}}]
  \label{def:appendix-1}
  A \define{Cartesian tangent category with scalar} $R$ \define{multiplication} consists of a Cartesian tangent category $\cat{C}$, a commutative ring object $R$, and an $R$-module structure on $T$ that is tangent-stable and for which some diagrams commute.
\end{definition}

For more details, such as the diagrams that need to commute, see \cite[Definition 2.24]{AinBloh25}. 

\subsection{Groups in tangent categories}
\label{sec:groups-tang-categ}

We now fix $\cat{C}$ to be a Cartesian tangent categories with scalar $R = \R$ multiplication. In particular, we require that $\R$ is an object of $\cat{C}$. Group objects $\cat{C}$ admit Lie algebras. This is a consequence of the much more general main result of \cite{AinBloh25}, which states that differentiable groupoid objects in Cartesian tangent categories with $R$ multiplication differentiate to a generalization of Lie algebroids.

\begin{theorem}[{\cite[Theorem 7.5]{AinBloh25}}]
  \label{thm:appendix-2}
Let $G$ be a group object in $\cat{C}$, a Cartesian tangent category with scalar $\R$-multiplication. Then the set of left-invariant vector fields $\fk{X}(G)^L$ identifies with the sections $\Gamma(*, \fk{g})$, where $\fk{g}$ is the pullback of $TG$ along the unit $* \to G$. The Lie bracket on $\fk{X}(G)$ induces a Lie bracket on $\fk{X}(G)^L$ that makes the $\R$-module $\fk{X}(G)^L$ a Lie algebra. 
\end{theorem}

\begin{remark}
  In this result, the $\R$-modules $\fk{X}(G)^L$ and $\Gamma(*, \fk{g})$ are not viewed as objects of $\cat{C}$. Indeed, in this generality there is no reason that they should be objects of $\cat{C}$.
\end{remark}

We will elucidate some of the ingredients and consequences of this theorem.  In what follows, we will identify $T(X\times Y) \cong TX \times TY$. Given $f\colon X \times Y \to Z$, we have the partial tangent morphisms
\begin{align*}
  &
    \begin{tikzcd}[ampersand replacement = \&]
      T_{(1)}f\colon TX \times Y \ar[r, "1 \times 0"] \& TX \times TY \ar[r, "Tf"] \& TZ 
    \end{tikzcd} \\
  &    \begin{tikzcd}[ampersand replacement = \&]
    T_{(2)}f\colon X \times TY \ar[r, "0 \times 1"] \& TX \times TY \ar[r, "Tf"] \& TZ.
  \end{tikzcd} \\  
\end{align*}
Thus $T_{(1)}T_{(2)}f$ and $T_{(2)}T_{(1)}f$ are two morphisms $TX \times TY \to T^2Z$. They are related by the canonical flip:
\begin{equation*}
  T_{(1)}T_{(2)}f = \tau_Z \circ  T_{(2)}T_{(1)}f.
\end{equation*}
Indeed, by naturality of $\tau$, and its compatibility with the Cartesian structure,
\begin{align*}
  &\tau_Z \circ T_{(2)}T_{(1)}f = \tau_Z \circ T^2f \circ (0, T0) = T^2f \circ (\tau_Z \circ 0, \tau_Z \circ T0) \\
  &\qquad = T^2f \circ (T0, 0) = T_{(1)}T_{(2)}f.
\end{align*}
We have the usual expression of the total derivative as the sum of partial derivatives.
\begin{lemma}[{\cite[Proposition 2.10]{CocCrut14}}]
  \label{lem:appendix-1}
  For $f\colon X \times Y \to Z$, we have
  \begin{equation*}
    Tf = T_{(1)}f \circ (1 \times \pi) + T_{(2)}f \circ (\pi \times 1).
  \end{equation*}
\end{lemma}

We will also need the following technical result.
\begin{lemma}
  \label{lem:appendix-3}
  Suppose that $v\colon X \to TX$ and $w\colon Y \to TY$ are vector fields, and $a\colon A \to X \times Y$ is a morphism. If
  \begin{equation*}
    Tf \circ (0,w) \circ a = Tf \circ (v,0) \circ a = 0 \circ f \circ a,
  \end{equation*}
  then
  \begin{equation*}
    T_{(1)}T_{(2)}f \circ (v,w) \circ a = T_{(2)}T_{(1)} f \circ (v,w) \circ a. 
  \end{equation*}
\end{lemma}
\begin{proof}
  Since $\lambda = \lambda \tau$, it suffices to show that $T_{(1)}T_{(2)}f \circ (v,w) \circ a = \lambda \circ \varepsilon$ for some map $\varepsilon\colon A \to TZ$. To do this, we use the fact that $\lambda$ is a kernel. Compute
  \begin{align*}
    \pi_{TZ} \circ T_{(1)}T_{(2)}f &= \pi_{TZ} \circ T(T_{(2)}f) \circ (1 \times 0) \\
                                               &= T_{(2)}f \circ (\pi_X \times \pi_{TY}) \circ (1 \times 0) \\
                                   &= T_{(2)}f \circ (\pi_X \times 1_{TY}) \\
                                   &= Tf \circ (0 \times 1) \circ (\pi_X \times 1_{TY}) \\
                                   &= Tf \circ (0 \circ \pi_X, 1_{TY}) \\
    T\pi_Z \circ T_{(1)}T_{(2)}f &= T(\pi_Z \circ T_{(2)}f) \circ (1 \times 0) \\
                                   &= Tf \circ (1 \times T\pi_Y) \circ (1 \times 0) \\
                                   &= Tf \circ (1_{TX} \times 0 \circ \pi_Y)
  \end{align*}
  These evaluate to $Tf \circ (0,w) \circ a$ and $Tf \circ (v,0) \circ a$ when pulled back along $(v,w) \circ a$, respectively. By assumption, these coincide with $0 \circ f \circ a$. Since $\lambda$ is a kernel, there must exist $\varepsilon\colon a \to TZ$ with $\lambda \circ \varepsilon = T_{(1)}T_{(2)}f \circ (v,w) \circ a$. Then
  \begin{equation*}
    T_{(2)}T_{(1)}f \circ (v,w) \circ a
    = \tau \circ T_{(1)}T_{(2)}f \circ (v,w) \circ a
    = \tau \circ \lambda \circ \varepsilon
    = \lambda \circ \varepsilon
    = T_{(1)}T_{(2)}f \circ (v,w) \circ a. 
  \end{equation*}
\end{proof}

Let $G$ be a group object in $\cat{C}$, with multiplication $m$, unit $e$, and inversion $i$. Then $TG$ is also a group object in $\cat{C}$, with multiplication $T m$, unit $Te = 0 \circ e$, and inversion $T i$.

\begin{definition}
  The \define{left-invariant} vector fields on $G$ are those $v \in \fk{X}(G)$ such that the following commutes:
  \begin{equation*}
    \begin{tikzcd}
      G \times TG \ar[r, "T_{(2)}m"] & TG \\
      G \times G \ar[u, "1 \times v"] \ar[r, "m"] & G. \ar[u, "v"] 
    \end{tikzcd}
  \end{equation*}
  We write $\fk{X}(G)^L$ for the collection left-invariant vector fields on $G$. 
\end{definition}
It takes some work to verify that the bracket of two left-invariant vector fields is again left-invariant. This is achieved in \cite[Theorem 6.6]{AinBloh25}, which applies more generally to groupoid objects. Let $\fk{g}$ denote the pullback of $TG$ along $e$. We have the left trivialization
\begin{equation*}
\varphi\colon G \times \fk{g} \hookrightarrow G \times TG \xar{T_{(2)}m} TG,
\end{equation*}
and its inverse, denoted $\psi$, determined by
\begin{equation*}
  \begin{tikzcd}
    TG  & G \times TG  & G \times TG \\
    & & G \times \fk{g}.
    \ar["{(i\circ \pi, 1)}", from=1-1, to=1-2]
    \ar["{1 \times T_{(2)}m}", from=1-2, to=1-3]
    \ar["\psi"', from=1-1, to=2-3]
    \ar[hook, from=2-3, to=1-3]
  \end{tikzcd}
\end{equation*}

By \cite[Theorem 7.4]{AinBloh25}, the $\R$-module of sections $\Gamma(*,\fk{g})$ is isomorphic to the $\R$-module $\fk{X}(G)^L$. The isomorphism and its inverse are respectively given by
  \begin{equation*}
\xi \mapsto \varphi \circ (1 \times \xi), \textrm{ and } v \mapsto v \circ u.
  \end{equation*}

\begin{remark}
  If $G$ is an elastic diffeological group, then $\fk{g} = T_eG$ carries the subspace diffeology. Furthermore, since $\cat{Dflg}$ is locally Cartesian closed, the space of left-invariant vector fields $\fk{X}(G)^L$ inherits the subspace diffeology of the functional diffeology on $\fk{X}(G)$. The isomorphism of $\R$-modules $\fk{g} \cong \fk{X}(G)^L$ is also an isomorphism in $\cat{Dflg}$. 
\end{remark}

There is another way to view the bracket on left-invariant vector fields, which we have not seen expressed in the literature. We let $c$ denote the conjugation
\begin{equation*}
  c\colon G \times G \to G, \quad c \coloneqq m \circ (m, i \circ \pr_1).
\end{equation*}
\begin{lemma}
  \label{lem:appendix-2}
  Let $v$ and $w$ be left-invariant vector fields. Then
  \begin{equation*}
    T_{(1)}T_{(2)}c \circ (v,w) \circ e = \delta(v,w) \circ e.
  \end{equation*}
\end{lemma}
This amounts to the classic expression of the Lie bracket as the second derivative of the conjugation. 
\begin{proof}
We break the computation into parts. First, by expanding definitions and using associativity of arrow composition,
  \begin{align*}
    T_{(2)}c &= Tm \circ (T_{(2)}m,\ T_{(2)}(i \circ \pr_1)) \\
             &= Tm \circ (T_{(2)}m, \ 0 \circ i \circ \pr_1) &&\text{using } T_{(2)}(i \circ \pr_1) = Ti \circ 0 \circ \pr_1 \\
    & &&\text{and } Ti \circ 0 = 0 \circ i \\
    &= T_{(1)}m \circ (T_{(2)}m, i \circ \pr_1).
  \end{align*}
Similarly,
  \begin{align*}
    T_{(1)}T_{(2)}c &= TT_{(1)}m \circ \Big(T_{(1)}T_{(2)}m, \ T_{(1)} (i \circ \pr_1)\Big) \\
    &= TT_{(1)}m \circ \Big(T_{(1)}T_{(2)}m, \ T(i \circ 0) \circ \pr_1 \Big).
  \end{align*}
Using Lemma \ref{lem:appendix-1}, we decompose the right side into a sum
  \begin{align*}
    T_{(1)}T_{(2)}c &= T_{(1)}T_{(1)}m \circ \Big(T_{(1)}T_{(2)}m,\ i \circ 0 \circ \pi \circ \pr_1 \Big) \\
    &\qquad + T_{(2)}T_{(1)}m \circ \Big(\pi \circ T_{(1)}T_{(2)}m, \ T(i \circ 0) \circ \pr_1\Big).
  \end{align*}
  Now we pre-compose with $(v,w) \circ e$. The most substantive term is
  \begin{align}
    T_{(1)}T_{(2)}m \circ (v,w) \circ e &= TT_{(2)}m \circ (v, 0 \circ w) \circ e \nonumber \\
                                        &= TT_{(2)}m \circ (1 \times Tw) \circ (v, 0) \circ e \nonumber \\
                                        &= Tw \circ Tm \circ (v \circ e, 0 \circ e) \label{eq:appendix-1}\\
                                        &= Tw \circ v \circ e.  \label{eq:appendix-2}
  \end{align}
  For \eqref{eq:appendix-1}, apply $T$ to the diagram expressing left-invariance of $w$; for \eqref{eq:appendix-2}, use the fact that $0 \circ e$ is the unit of $TG$, and $Tm$ is the multiplication. Therefore
  \begin{align*}
    T_{(1)}T_{(2)}c \circ (v,w) \circ e &= T_{(1)}T_{(1)}m \circ (Tw \circ v \circ e, e) \\
    &\qquad + T_{(2)}T_{(1)}m \circ (w \circ e, Ti \circ v \circ e).
  \end{align*}
  The first summand expands to
  \begin{equation*}
    T_{(1)}T_{(1)}m \circ (Tw \circ v \circ e, e) = T^2m(Tw \circ v \circ e, 0 \circ 0 \circ e),
  \end{equation*}
  and $0 \circ 0 \circ e$ is the group unit in $T^2G$, we get $Tw \circ v \circ e$.  For the second summand, recall that $\tau \circ T_{2}T_{(1)} = T_{(1)}T_{(2)}$. We also have
  \begin{equation*}
    Ti \circ v \circ e = -v \circ e,
  \end{equation*}
  because $Ti$ is the group inversion on $TG$, and $Tm \circ (a \circ e, b \circ e) = a \circ e + b \circ e$. Putting these together, the second summand is
  \begin{equation*}
    T_{(2)}T_{(1)}m \circ (w \circ e, Ti \circ v \circ e) = - \tau \circ Tv \circ w \circ e.
  \end{equation*}
  We conclude that
  \begin{equation*}
    T_{(1)}T_{(2)}c \circ (v,w) \circ e = \delta(v,w) \circ e,
  \end{equation*}
  as desired.
\end{proof}

\printbibliography

\end{document}1